\documentclass[twoside]{amsart}
\usepackage{nuclear, format,whiteNoise} 

 \renewcommand{\include}{\input}
 
\begin{document}

\title{Countably-Normed Spaces, Their Dual, and the Gaussian Measure}

\author{  Jeremy J. Becnel }
  \address{
  Department of Mathematics\\
  Louisiana State University\\
  Baton Rouge, LA 70803, USA\\
  e-mail: \sl  becnel@gmail.com}
\date{\today}


\begin{abstract} 
Here we present an overview of \CNS s.  We discuss the main topologies---weak, strong, and inductive---placed on the dual of a \CNS\ and discuss the \sfield s generated by these topologies. In particular, we show that under certain conditions the strong and inductive topologies coincide and the \sfield s generated by the weak, strong, and inductive topologies are equal. With this \sfield, we develop a Gaussian measure on the dual of a nuclear space. The purpose in mind is to provide the background material for many of the results used is White Noise Analysis.
\end{abstract}

\maketitle

\section{Topological Vector Spaces} \label{S:tvs}
In this section we review the basic notions of topological vector spaces along and provide proofs a few useful results.

\subsection{Topological Preliminaries} \label{SS:toppre}
Let $E$ be a real vector space.

A \textit{vector topology} $\tau$ on $E$ is a   topology such that
addition  $E\times E\to E:(x,y)\mapsto x+y$  and scalar
multiplication $\R \times E\to E: (t,x)\mapsto tx$ are
continuous. If $E$ is a complex vector space we require that
$\C\times E\to E: (\alpha, x)\mapsto\alpha x$ be
continuous.

It is useful to observe that when $E$ is equipped with a vector
topology, the  translation maps
$$t_x: E\to E:y\mapsto y+x$$
are continuous, for every $x\in E$, and are hence also
homeomorphisms since $t_x^{-1}=t_{-x}$.

A \textit{topological vector space} is a vector space equipped
with a vector topology.

Recall that a \textit{local base} of a vector topology $\tau$ is a family of
open sets $\{U_{\alpha}\}_{\alpha\in I}$ containing $0$ such that
if $W$ is any open set containing $0$ then $W$ contains some
$U_{\alpha}$. A set $W$ that contains an open set containing $x$ is called a \emph{neighborhood} of $x$. If $U$ is any open set and $x$ any point in $U$ then
$U-x$ is an open neighborhood of $0$ and hence contains some
$U_\alpha$, and so $U$ itself contains a neighborhood
$x+U_{\alpha}$ of $x$: 
\begin{equation}\label{UUalpha}
\hbox{\textit{If $U$ is open and $x\in U$ then
$x+U_{\alpha}\subset U$, for some $\alpha\in I$.}}
\end{equation} 
Doing this for each point $x$ of $U$, we see that
each open set is the union of translates of the local base sets
$U_{\alpha}$.

If $\mathcal{U}_x$ denotes the set of all neighborhoods of a point $x$ in a topological  space $X$, then $\mathcal{U}_x$ has the following properties:
\begin{numList}
  \item $x \in U$ for all $U \in \mathcal{U}_x$.
  \item if $U \in \mathcal{U}_x$ and $V \in \mathcal{U}_x$, then $U \cap V \in  \mathcal{U}_x$.
  \item if $U \in \mathcal{U}_x$ and  $U \subset V$, then $V \in \mathcal{U}_x$.
  \item if $U \in \mathcal{U}_x$, then there is some $V \in \mathcal{U}_x$ with $U \in    \mathcal{U}_y$ for all $y \in V$. (taking $V$ to be the interior of $U$ is sufficient).
\end{numList}
Conversely if $X$ is any set and a non-empty collection of subsets $\mathcal{U}_x$ is given for each $x \in X$, then when the conditions above are satisfied by the $\mathcal{U}_x$, exactly one topology can be defined on $X$ in such a way to make  $\mathcal{U}_x$ the set of neighborhoods of $x$ for each $x\in X$. A set $V \subset X$ is called \emph{open} if for each $x \in V$, there is a $U \in \mathcal{U}_x$ with $U \subset V$. \cite{aR64} 

In most cases of interest a topological vector space has a local base consisting of convex sets. We call such spaces \emph{locally convex topological vector spaces}. 

In a topological vector space there is the notion of bounded sets. A set $D$ in a topological vector space is said to be \emph{bounded}, if for every neighborhood $U$ of $0$ there is some $\lambda > 0$ such that $D \subset \lambda U$. If $\{U_{\alpha}\}_{\alpha\in I}$ is a local base, then it is easily seen that $D$ is bounded if and only if to each $U_{\alpha}$ there corresponds $\lambda_{\alpha} > 0$ with $D \subset \lambda_{\alpha} U_{\alpha}$. \cite{aR64}

A set $A$ in a vector space $E$ is said to be \emph{absorbing} if given any $x \in E$ there is an $\eta$ such that $x \in \l A $ for all $ |\l| \geq \eta$. The set $A$ is called \emph{balanced} if, for all $ x\in A$, $\l x \in A$ whenever $| \l | \leq 1$. Also, a set $A$ in a vector space $E$ is call \emph{symmetric} if $-A = A$. Finally, although the next concept is very common, the term we use for it is not, so we make a formal definition:

\begin{Def} \label{D:lpc}
A subset $A$ of a topological space $X$ is \emph{limit point compact} if every infinite subset of $A$ has a limit point.
\end{Def}

\begin{Remark}
The term \emph{limit point compact} is not the standard term for spaces with the above property. In fact, I do not believe there is a standard term. I have seen it called ``Fr\'echet compactness", ``relative sequential compactness", and the ``Bolzano-Weierstrass property". The term \emph{limit point compact} was taken directly from Munkres \cite{jM00}. It is my personal favorite term; at the very least it is descriptive.
\end{Remark}

\subsection{Bases in Topological Vector Spaces}
Here we take the time to prove some general, but very useful, results about local bases for topological vector spaces. Most of the results in this subsection are taken from Robertson \cite{aR64}.

\begin{Lemma} \label{L:balance}
Every topological vector space $E$ has a base of balanced neighborhoods.
\end{Lemma}
\begin{proof}
Let $U$ be a neighborhood of $0$ in $E$. Consider the function $h:\C \times E \to E$ given by $h(\l,x) = \l x$. Since $E$ is a topological vector space, $h$ is continuous at $\l =0, x = 0$. So there is a neighborhood $V$ and $\e > 0$ with $\l x \in U$ for $|\l| \leq \e$ and $x \in V$. Hence $\l V \subset U$ for $|\l| \leq \e$. Therefore $ \fr(\e/\alpha) V \subset U$ for all $\alpha$ with $| \alpha | \geq 1$. Thus  $\e V \subset U' = \bigcap_{|\alpha| \geq 1} \alpha U \subset U$. Now since $V$ is a neighborhood of $0$ so is $\e V$. Hence $U'$ is a neighborhood of $0$. If $x \in U'$ and $0 \leq |\l| \leq 1$, then for $|\alpha| \geq 1$, we have $x \in \fr(\alpha/\l) U$ (since $\left| \fr(\alpha/\l) \right| \geq 1$). So $\l x \in \alpha U$ for $| \l | \leq 1$. Hence $\l x \in U'$. Therefore $U'$ is balanced.
\end{proof}

\begin{Lemma} \label{L:convexBase}
Let $E$ be a vector space. Let $\mathcal{B}$ be a collection of subsets of $E$ satisfying:
\begin{romanList}
 \item if $U,V \in \mathcal{B}$, then there exist $W \in \mathcal{B}$ with $W \subset U\cap V$.
 \item if $U \in \mathcal{B}$ and $\l \neq 0$, then $ \l U \in \mathcal{B}$.
 \item if $U \in \mathcal{B}$, then $U$ is balanced, convex, and absorbing.
\end{romanList}
Then there is a topology making $E$ a locally convex topological vector space with $\mathcal{B}$ the base of neighborhoods of $0$.
\end{Lemma}

\begin{proof}
Let $\A$ be the set of all subsets of $E$ that contain a set of $\mathcal{B}$. For each $x$ take $x + \A$ to be the set of neighborhoods of $x$. We need to see that (1)-(4) are satisfied from subsection \ref{SS:toppre}. 

For (1), we have to show $x \in A$ for all $ A \in x + \A$. Note that since each $U \in \mathcal{B}$ is absorbent, there exists a non-zero $\l$ such that $0 \in \l U$. But then $0 \in \l^{-1} \l U = U$. So each $U \in \mathcal{B}$ contains $0$. So $x \in A$ for all $A \in x + \A$.

For (2), we have to show that if $A,B \in \A$, then $(x + A)\cap (x+ B) \in x + \A$ for each $x \in E$. Recall that $U \subset A$ and $ V \subset B$ for some $U,V \in \mathcal{B}$. So $U \cap V \subset A \cap B$. By the first hypothesis, there is a $W \in \mathcal{B}$ with $W \subset U \cap V \subset A \cap B$. Thus $A \cap B \in \A$ and hence $(x + A)\cap (x+ B) \in x + \A$ for each $x \in E$.

Next (3) is clear from the definition of $\A$, since if $A \in \A$ and $A \subset B$ then $B \in  \A$. 

Finally for (4), we must show that if $x + A \in x + \A$, then there is an $V \in x + \A$ with $x+ A \in y + \A$ for all $y \in V$. If $A \in \A$ take a $U \in \mathcal{B}$ with $U \subset A$. Now we see that $x + A$ is a neighborhood of each point $y \in  x + \fr(1/2) U $.  Since $ y \in x + \fr(1/2) U$ we have $y-x \in \fr(1/2) U$. Thus $y-x + \fr(1/2) U \subset \fr(1/2) U + \fr(1/2) U \subset A$. Hence $y + \fr(1/2) U \subset x+ A$. Thus $x -y + A \supset \fr(1/2) U$. So $x-y+A \in \A$. Therefore $ y + x - y + A = x + A \in y + \A$.

To prove continuity of addition, let $U \in \mathcal{B}$. Then if $x \in a + \fr(1/2) U$ and $ y \in b + \fr(1/2) U$, we have $x + y \in  a + b + U$.

Finally, to see that scalar multiplication, $\l x$,  is continuous at $ x =a,  \l = \alpha$, we should find $\delta_1$ and $ \delta_2 $ such that $\l x - \alpha a \in U$ whenever $| \l - \alpha | < \delta_1 $ and $x \in a + \delta_2 U$. Since $U$ is absorbing, there is a $\eta$ with $a \in \eta U $. Take $\delta_1$ so that $ 0 < \delta_1 < \fr(1/{2\eta})$ and take $\delta_2$ so that
 $0< \delta_2 < \fr(1/{2 \left( |\alpha| + \delta_1 \right)})$. Now observe
 \begin{align*}
 \l x - \alpha a &= \l (x -a) + (\l - \alpha) a \\
  & \in (|\alpha| + \delta_1) \delta_2 U + \delta_1 \eta U \\
  & \subset \fr(1/2) U + \fr(1/2) U \subset U.  
\end{align*}
Thus we are done.
\end{proof}

\subsection{Topologies generated by families of topologies} \label{SS:famtop}

Let $\{\tau_{\alpha}\}_{\alpha\in I}$ be a collection of
topologies on a space.  It is natural and useful to consider
the least upper bound topology $\tau$, i.e. the coarsest topology
containing all sets of $\cup_{\alpha\in I}\tau_{\alpha}$. In our
setting, we work with each $\tau_{\alpha}$ a vector topology on a
vector space $E$.
\begin{Thm} \label{T:topfam}
The least upper bound topology $\tau$ of
a collection $\{\tau_{\alpha}\}_{\alpha\in I}$  of vector
topologies is again a vector topology. If $\{W_{\alpha, i}\}_{i\in
I_{\alpha}}$ is a local base for $\tau_{\alpha}$ then a local base
for $\tau$ is obtained by taking all finite intersections of the
form $W_{\alpha_1,i_1}\cap\cdots\cap W_{\alpha_n,i_n}$.
\end{Thm}

\begin{proof} Let $\mathcal{B}$ be the collection of all sets
which are of the form $W_{\alpha_1,i_1}\cap\cdots\cap
W_{\alpha_n,i_n}$.

Let $\tau'$ be the collection of all sets which are unions of
translates of sets in $\mathcal{B}$ (including the empty union).
Our first objective is to show that $\tau'$ is a topology on $E$.
It is clear that $\tau'$ is closed under unions and contains the
empty set. We have to show that the intersection of two sets in
$\tau'$ is in $\tau'$. To this end, it will suffice to prove the
following:
\begin{eqnarray}
\hbox{If $C_1$ and $C_2$ are sets in $\mathcal B$, and $x$ is a
point
in}  \nonumber\\
\hbox{the intersection of the translates $a+C_1$ and
$b+C_2$,}\label{toshow}\\
\hbox{ then $x+C\subset (a+C_1)\cap (b+C_2)$ for some $C$ in
$\mathcal B$.}\nonumber
\end{eqnarray}

 Clearly, it suffices to consider finitely many topologies
$\tau_{\alpha}$.
 Thus, consider vector topologies $\tau_1,...,\tau_n$ on $E$.

   Let $\mathcal{B}_n$ be the collection of
 all sets of the form $B_1\cap\cdots\cap B_n$ with
 $B_i$ in a local base for $\tau_i$, for each
 $i\in\{1,...,n\}$. We can check that if $D, D'\in\mathcal{B}_n$ then
 there is an $G\in \mathcal{B}_n$ with $G\subset D\cap D'$.

 Working with $B_i$ drawn from a given local base for $\tau_i$, let
  $z$ be a point in the intersection $B_1\cap \cdots\cap B_n$. Then
 there exist sets $B'_i$, with each $B'_i$ being in the local base
 for $\tau_i$, such that $z+B'_i\subset B_i$ (this follows from
 our earlier observation (\ref{UUalpha})). Consequently,
 $$z+\cap_{i=1}^nB'_i\subset  \cap_{i=1}^nB_i.$$
 Now consider sets $C_1$ and $C_2$, both in $\mathcal{B}_n$.
 Consider $a, b\in E$ and suppose $x\in (a+C_1)\cap (b+C_2)$.
 Then since $x-a\in C_1$ there is a set $C'_1\in \mathcal{B}_n $
 with $x-a+C'_1\subset C_1$; similarly, there is a
$C'_2\in \mathcal{B}_n $
 with $x-b+C'_2\subset C_2$. So $x+C'_1\subset a+C_1$
 and $x+C'_2\subset b+C_2$. So
 $$x+C\subset (a+C_1)\cap (b+C_2),$$
 where $C\in\mathcal{B}_n$ satisfies $C\subset C_1\cap C_2$.

 This establishes (\ref{toshow}), and shows that the intersection of
  two sets in $\tau'$ is in
 $\tau'$.

 Thus $\tau'$ is a topology. The definition of   $\tau'$ makes it
 clear that $\tau'$ contains each $\tau_{\alpha}$. Furthermore,
 if any topology $\sigma$ contains each $\tau_{\alpha}$ then all the
sets
 of $\tau'$ are also open relative to $\sigma$. Thus
 $$\tau'=\tau,$$
 the topology generated by the topologies $\tau_{\alpha}$.

 Observe that we have shown that if $W\in\tau$ contains $0$ then
 $W\supset B$ for some $B\in \mathcal{B}$.

 Next we have to show that $\tau$ is a vector topology. The
 definition of $\tau$ shows that $\tau$ is translation invariant,
 i.e. translations are homeomorphisms. So, for addition, it will
 suffice to show that addition $E\times E \to E: (x,y)\mapsto
 x+y$ is continuous at $(0,0)$.  Let $W\in\tau$ contain $0$.
 Then there is a $B\in\mathcal{B}$ with $0\in B\subset W$.
 Suppose $B=B_1\cap\cdots\cap B_n$, where each $B_i$ is in
 the given local base for $\tau_i$. Since $\tau_i$ is a vector
 topology, there are open sets $D_i, D'_i\in\tau_i$, both
 containing $0$, with
 $$D_i+D'_i\subset B_i.$$
  Then choose $C_i, C'_i$ in the local base for
 $\tau_i$ with $C_i\subset D_i$ and $C'_i\subset D'_i$.  Then
 $$C_i+C'_i\subset B_i.$$
 Now let $C=C_1\cap\cdots\cap C_n$, and $C'=C'_1\cap\cdots\cap
 C'_n$. Then $C, C'\in\mathcal{B}$ and $C+C'\subset B$. Thus,
 addition is continuous at $(0,0)$.

 Now consider the multiplication map $\R \times E\to E:
 (t,x)\mapsto tx$.   Let $(s,y), (t,x)\in \R \times E$. Then
 $$sy-tx= (s-t)x+t(y-x)+ (s-t)(y-x).$$
Suppose $F\in\tau$ contains $tx$. Then
$$F\supset tx+W',$$ for some
 $W'\in\mathcal{B}$. Using continuity of the addition
 map
 $$E\times E\times E\to E: (a,b,c)\mapsto a+b+c$$
 at $(0,0,0)$,
 we can choose $W_1, W_2, W_3\in\mathcal{B}$ with
 $W_1+W_2+W_3\subset W'$. Then we can choose $W\in \mathcal{B}$, such
that
  $$W\subset W_1\cap W_2\cap W_3$$ Then
 $W\in\mathcal{B}$ and
 $$W+W+W\subset W'.$$
 Suppose $W=B_1\cap\cdots\cap B_n$, where each $B_i$ is
 in the given local base for the vector topology $\tau_i$. Then for $s$
close enough to $t$,
 we have $(s-t)x\in B_i$ for each $i$, and hence $(s-t)x\in W$.
 Similarly, if $y$ is $\tau$--close enough to $x$ then $t(y-x)\in W$.
 Lastly, if $s-t$ is close enough to $0$ and $y$ is close enough
 to $x$ then $(s-t)(y-x)\in W$. So $sy-tx\in W'$, and so
 $sy\in F$, when $s$ is close enough to $t$ and $y$ is
 $\tau$--close enough to $x$.
\end{proof}

The above result makes it clear that if each $\tau_\alpha$ has a
convex local base then so does $\tau$. Note also that if at least
one $\tau_\alpha$ is Hausdorff then so is $\tau$.

A family of topologies $\{\tau_{\alpha}\}_{\alpha\in I}$ is
\textit{directed} if for any $\alpha, \beta\in I$ there is a
$\gamma\in I$ such that
$$\tau_{\alpha}\cup\tau_{\beta}\subset\tau_{\gamma}.$$
In this case every open neighborhood of $0$ in the generated
topology contains an open neighborhood in one of the topologies
$\tau_{\gamma}$.


\section{Countably--Normed Spaces}
We begin with the basic definition of a \emph{\CNS} and a \emph{\CHS}.

\begin{Def}
Let $V$ be a topological vector space over \C\ with topology given by a family of norms $\{\norm|\cdot|_n ; n = 1,2,\ldots \}$. Then $V$ is a \emph{\CNS}. The space $V$ is called a \emph{\CHS} if each \norm|\cdot|_n is an inner product norm and $V$ is complete with respect to its topology.
\end{Def}

\begin{Remark} \label{R:inc}
By considering the new norms $\Vert v \Vert_n = (\sum_{k=1}^n {\norm|v|_k}^2)^{\fr(1/2)}$ we may assume that the family of norms
$\{\norm|\cdot|_n ; n = 1,2,\ldots \}$ is increasing, i.e.
\[
\norm|v|_1 \leq \norm|v|_2 \leq \cdots \leq \norm|v|_n \leq \cdots, 
\forall v \in V
\]
\end{Remark}

If $V$ is a \CNS, we denote the completion of $V$ in the norm \norm|\cdot|_n by $V_n$. Then $V_n$ is by definition a Banach space. Also in light of Remark~\ref{R:inc} we can assume that 
\[
V \subset \cdots \subset V_{n+1} \subset V_{n} \subset \cdots \subset V_1
\]

\begin{Lemma} \label{L:iCont}
The inclusion map from $V_{n+1}$ into $V_n$ is continuous. 
\end{Lemma}

\begin{proof}
Consider an open neighborhood of 0 in $V_n$ given by
\[
B_n(0,\epsilon) = \{ v \in V_n ; \norm|v|_n < \epsilon \}
\]
Let $i_{n+1,n}:V_{n+1} \to V_{n}$ be the inclusion map. Now
\[
i_{n+1,n}^{-1}(B_n(0,\epsilon)) = 
\{ v \in V_{n+1} ; \norm|v|_n < \epsilon \} 
\supset B_{n+1}(0, \epsilon) \text{ since } \norm|v|_n  \leq \norm|v|_{n+1}
\]
Therefore $i_{n+1,n}$ is continuous.
\end{proof}

\begin{Prop} \label{P:complete}
Let $V$ be a \CNS. Then $V$ is complete if and only if $V = \bigcap_{n=1}^{\infty} V_n$.
\end{Prop}

\begin{proof}
Suppose $V = \bigcap_{n=1}^{\infty} V_n$ and $\{v_k\}_{k=1}^{\infty}$ is Cauchy in $V$. By definition
$\{v_k\}_{k=1}^{\infty}$ is Cauchy in $V_n$ for all $n$. Since $V_n$ is complete, a limit $v^{(n)}$ exist in $V_n$. 
Using that the inclusion map $i_{n+1,n}:V_{n+1} \to V_{n}$ is continuous (by \Lemmaref{L:iCont}) and that
\[
V \subset \cdots \subset V_{n+1} \subset V_{n} \subset \cdots \subset V_1
\]
we have that all the $v^{(n)}$ are the same and belong to each $V_n$. 
Thus they are in $V = \bigcap_{n=1}^{\infty} V_n$.
Let us call this element $v \in V$.

Since \norm|v_k-v^{(n)}|_m $\to 0$ for all $m$ we have that \norm|v_k-v|_m $\to 0$ for all $m$. 
Hence $ v = \lim_{k \to \infty} v_k$ in $V$. Thus $V$ is complete.

Conversely, let $V$ be complete and take $v \in \bigcap_{n=1}^{\infty} V_n$. We need to show $v$ is in $V$. For each $n$ we can find $v_n \in V$ such that $\norm|v - v_n|_n < \fr(1/n)$ (using that $V$ is dense in $V_n$). Now for any $k < n$ we have $\norm|v-v_n|_k \leq \norm|v-v_n|_n < \fr(1/n)$. Thus $\lim_{n \to \infty} \norm|v-v_n|_k = 0$. This gives us that $\{v_n\}$ is Cauchy with respect to all norms $\norm|\cdot|_k$ where $k = 1,2,\ldots$ 

Let $\overline{v} = \lim_{n \to \infty} v_n$ in $V$. Since for all $k$ we have $\overline{v}, v \in V_k$ and $\lim_{n \to \infty} \norm|\overline{v} - v_n|_k = 0$, we see that $v=\overline{v}$. Thus $v \in V$ and we have $V \supset \bigcap_{n=1}^{\infty} V_n$. That $V \subset \bigcap_{n=1}^{\infty} V_n$ is obvious, since $V \subset V_n$ for all $n$.
\end{proof}

\subsection{Open Sets in $V$} \label{SS:openV}
In light of \Thmref{T:topfam}, we see that a local base for $V$ is given by sets of the form:
\[
 B = \B_{n_1}( \e_1 ) \cap \B_{n_2}(\e_2) \cap \cdots \cap \B_{n_k}(\e_k),
\]
where $\B_{n_i}(\e_i)= \left\{ v \in V ; \norm|v|_{n_i} < \e_i  \right\}$ is the $\norm|\cdot|_{n_i}$--ball of radius $\e_i$ in $V$. 

\begin{Prop} \label{P:base}
Let $V$ be a \CNS. For every element $B$ of the local base for $V$ there exist $n$ and $\e > 0$ such that $\B_{n}(\e) \subset B$.
\end{Prop}

\begin{proof}
Let $B = \B_{n_1}( \e_1 ) \cap \B_{n_2}(\e_2) \cap \cdots \cap \B_{n_k}(\e_k)$ be an element of the local base for $V$. 
Then take 
$n = \max_{1 \leq j \leq k} n_j$ and 
$\e = \min_{1 \leq j \leq k} \e_j$.
Observe $\B_{n}(\e) \subset B$ since for $v \in \B_n(\e)$ we have 
$\norm|v|_{n_j} \leq \norm|v|_{n} < \e \leq \e_j$ for any $j \in \{1,2, \ldots, k\}$. Thus $v \in B$.
\end{proof}

\begin{Cor} \label{C:base}
Let $V$ be a \CNS. Then a local base for $V$ is given by the collection $\{ \B_{n}({\fr(1/k)}) \}_{n,k = 1}^{\infty}$.
\end{Cor}

\begin{Cor} \label{C:baseCount}
Let $V$ be a \CNS. Then a local base for $V$ is given by the collection 
$\{ \B_{k}({\fr(1/k)}) \}_{k = 1}^{\infty}$. 
Moreover we have that $\B_{1}(1) \supset  \B_{2}({\fr(1/2)}) \supset \cdots $
\end{Cor}

\begin{proof}
Let $U$ be a neighborhood of $0$. By \Corref{C:base} there are positive integers $n$ and $k$ such that 
$B_{n}({\fr(1/k)}) \subset U$. If $n \geq k$, we have that $\B_n({\fr(1/n)}) \subset \B_n({\fr(1/k)})$ since $\fr(1/n) \leq \fr(1/k)$. If $n \leq k$, then $\B_k({\fr(1/k)}) \subset \B_n({\fr(1/k)})$ 
since $\norm|v|_k < \fr(1/k)$ gives us that $\norm|v|_n \leq \norm|v|_k < \fr(1/k)$.

For $m \geq k$ we have that $\B_m({\fr(1/m)}) \subset \B_k({\fr(1/k)})$ since $\norm|v|_k \leq \norm|v|_m$ and
$\fr(1/m) < \fr(1/k)$.  
\end{proof}

\subsection{Bounded Sets in $V$} \label{SS:bndV}
Recall that a subset $D$ of a \CNS\ $V$ is said to be \emph{bounded} if for any neighborhood $U$ of zero in $V$ there is a positive number $\lambda$ such that $D \subset \lambda U$ (see subsection \ref{SS:toppre}). This leads us to the following useful proposition:

\begin{Prop} \label{P:bndIffV}
A set $D$ in a \CNS\ $V$ is bounded if and only if $\sup_{v \in D} \norm|v|_n < \infty$ for all $n \in \{ 1,2,\ldots \}$.
\end{Prop}

\begin{proof} ($\Rightarrow$) Suppose $D$ is a bounded set in $V$. Take the open neighborhood $\rB_{n}(1)= \left\{ v \in V ; \norm|v|_{n} < 1  \right\}$ in $V$. Since $D$ is bounded in $V$ there is an $\l > 0$  such that $D \subset \l \B_n(1)$. Thus $\sup_{v \in D} \norm|v|_n \leq \l$.

($\Leftarrow$) Suppose $U$ is a neighborhood of $0$ in $V$. Then by \Propref{P:base} there is an $\B_n(\e) \subset U$. Let $ \sup_{v \in D} \norm|v|_n = M < \infty $. Then $D \subset \fr(M+1/{\e}) \B_n(\e) \subset \fr(M+1/\e) U$. So $D$ is bounded.
\end{proof}

\subsection{The Dual} \label{SS:dualV}
Again take $V$ to be a \CNS\ associated with an increasing sequence of norms $\{\norm|\cdot|_n \}_{n=1}^{\infty}$ and let $V_n$ be the completion of $V$ with respect to the norm $\norm|\cdot|_n$. We denote the dual space of $V$ by $V'$. Let \ip<\cdot \, ,\cdot> denote the bilinear pairing of $V'$ and $V$.

Of course, each Banach space $V_n$ also has a dual, which we denote by $V_n'$. We use the notation to $\norm|\cdot|_{-n}$ to denote the operator norm on the Banach space $V_n'$. The relationship between $V'$ and each $V_n'$ is discussed in the next proposition.

\begin{Prop} \label{P:dualInc}
The dual of a \CNS\ $V$ is given by $V'\negthickspace=\negmedspace\bigcup_{n=1}^{\infty} V_n'$ and we have the inclusions 
\[ 
V_1' \subset \cdots \subset V_n' \subset V_{n+1}' \subset \cdots V'
\]
Moreover, for $f \in V_n'$ we have $\norm|f|_{-n} \geq  \norm|f|_{-n-1}$.
\end{Prop}

\begin{proof} ($\supset$) Take $v' \in V_n'$. Then $v'$ is continuous on $V_n$ with topology coming from the norm \norm|\cdot|_n . Thus $v'$ is continuous on $V$, since $V \subset V_n$ and the norm \norm|\cdot|_n is one of the norms generating the topology on $V$.

($\subset$) Take $v' \in V'$. Since $v'$ is continuous on $V$ the set 
\[
{v'}^{-1}\left(-1, 1\right) = \{v \in V \medspace ; \medspace \vert \ip<v',v> \vert < 1 \}
\]
is open in $V$. So we can find a member $B$ of the local base for $V$ such that $B \subset {v'}^{-1}\left(-1, 1\right)$. By to \Propref{P:base} we have that $\rB_n(\e) \subset {v'}^{-1}\left(-1, 1\right)$ for some positive integer $n$ and some $\e > 0$.

Thus for all $v \in V$ with $\norm|v|_n < \e$ we have that 
$| \ip<v',v>| < 1$. Since $V$ is dense in $V_n$,  if $v \in V_n$ and $\norm|v|_n \leq \e$ then $| \ip<v',v>| \leq 1$. Thus $v' \in V_n'$.

To see that $V_n' \subset V_{n+1}'$ take $f \in V_n'$. Then for all $v \in V_n$ we have that 
\[
\vert f(v) \vert \leq \norm|f|_{-n} \norm|v|_{n}  \leq \norm|f|_{-n} \norm|v|_{n+1}
\]
Since $V_{n+1} \subset V_n$, the above holds for all $v \in V_{n+1}$. Thus $f \in V_{n+1}'$ and $\norm|f|_{-n-1} \leq \norm|f|_{-n}$.
\end{proof}

\begin{Prop} \label{P:borno}
A linear functional $f$ on $V$ is continuous if and only if $f$ is bounded on bounded sets of $V$.
\end{Prop}

\begin{proof}
($\Rightarrow$) Let $f$ be a continuous linear functional on $V$. Then $f$ is in $V'$. So $f = \ip<v', \cdot>$ for some $v' \in V'$. Now by \Propref{P:dualInc}, $v' \in V_n'$ for some $n$.  Let $D \subset V$ be bounded. By \Propref{P:bndIffV} we have that $\sup_{v \in D} \norm|v|_n = M < \infty$. Using this we see that
$\sup_{v \in D} \aip<v',v> \leq M\norm|v'|_{-n}  < \infty$. Thus $f = \ip<v',\cdot>$ is bounded on bounded sets.

($\Leftarrow$) Suppose $f$ is bounded on bounded sets. Consider the local base sets 
$\B_1(1) \supset \B_2({\fr(1/2)}) \supset \cdots $ in $V$ as in \Corref{C:baseCount}. By contradiction we assume that $f$ is not in $V'$. Then $f$ is not in $V_k'$ for any $k$. So $f$ is not continuous on $V_k$ and hence not bounded on $\B_k( {\fr(1/k)} )$. Hence we can find a $v_k$ in $\B_k( {\fr(1/k)})$ such that  $|f(v_k)| > k$. The sequence $\{v_k\}_{k =1}^{\infty}$ goes to $0$ in $V$. Thus $\{v_k\}_{k =1}^{\infty}$ must be bounded. But then by hypothesis, 
$\{ f \! \left( v_k \right) \}_{k=1}^{\infty} $ should be bounded. But by construction it is not, a contradiction.
\end{proof}

\begin{Cor} \label{C:borno}
A linear functional $f$ on $V$ is continuous if and only if $f$ is bounded on some neighborhood of $0$ in $V$.
\end{Cor}

\begin{proof}
Suppose $f$ is bounded on some neighborhood $U$ of $0$ in $V$. Then for any $\alpha > 0$, $f$ is bounded on $\alpha U$. Let $D$ be a bounded set in $V$. Then $D \subset \l U$ for some $\l > 0$. So $f$ is bounded on $D$ and hence continuous by \Propref{P:borno}
\end{proof}

There are several topologies one can put on the dual space $V'$. The three most common are the weak, strong, and inductive topologies. In the following sections we discuss the properties of these three topologies and compare them against one another. Throughout this discussion, the topology on $V_n'$ is taken to be the usual strong topology (i.e. the topology induced by the operator norm on $V_n'$ as the dual of the Banach space $V_n$). 

\subsection{Bounded Sets of $V$ revisited} \label{SS:bndR}
Let $V$ be a \CNS. With the notion of the dual $V'$ of $V$ behind us (see subsection \ref{SS:dualV}), we can formulate a better understanding of bounded sets in $V$. We begin with the following simple definition:

\begin{Def} \label{D:weakBnd}
A set $D \subset V$ is said to be \emph{weakly bounded} if given a set 
$\N(v';\e)=\{ v\in V ; \aip<v',v> < \e \}$ there is a $\l > 0$ such that
$D \subset \l \N(v';\e)$.
\end{Def}

\begin{Thm} \label{T:bndEqV}
Suppose $V$ is a \CNS\ with dual $V'$. Let $D \subset V$. Then the following are equivalent:
\begin{numList}
 \item $D$ is bounded.
 \item $D$ is weakly bounded.
 \item The values of each $v' \in V'$ are bounded on $D$.
 \item For all $n$, we have $\sup_{v \in D} \norm|v|_n < \infty$.
\end{numList}
\end{Thm}

\begin{proof}
We have already shown that (1) and (4) are equivalent in \Propref{P:bndIffV}.

$( (1) \Rightarrow (2))$ Suppose $D$ is bounded in $V$. Take a $v' \in V'$. Then $v' \in V_n'$ for some $n$. For $v \in D$ we have $\aip<v',v> \leq \norm|v'|_{-n} \norm|v|_n \leq \norm|v'|_{-n} M_n$ where $M_n = \sup_{v \in D} \norm|v|_n$. Thus
we have $D \subset \fr({2 \norm|v'|_{-n} M_n}/ \e) \N(v';\e)$. So $D$ is weakly bounded.

$( (2) \Rightarrow (3))$
Suppose $D$ is weakly bounded in $V$. Take $v' \in V'$. By assumption $D \subset \l \N(v';\e)$ for some $\l > 0 $. So for $v \in D$ we have $\aip<v',v> \leq \l \e$.

$( (3) \Rightarrow (4))$
Consider $D \subset V \subset V_n$. By hypothesis all $v' \in V'$ are bounded on $D$.
In particular all $v' \in V_n' \subset V'$ are bounded on $D$. This means the linear functionals $\{ \ip<\cdot,v> ; v \in D \}$ are pointwise bounded on $V_n'$. Thus we can apply the uniform boundedness principle to see that $\sup_{v \in D} \norm|v|_n < \infty$.
\end{proof}

\subsection{The Metric on $V$} Let $V$ be a \CNS. Define the function $\p : V \times V \to [0, \infty)$ by 
\begin{equation} \label{E:metric}
\p(v,u) = \sum_{n=1}^{\infty} \fr(1/2^n) \fr({\norm|v-u|_n}/{1+\norm|v-u|_n}).
\end{equation}
First observe that $\p$ is a metric on $V$. From the above definition it is obvious that $\p(v,v) = 0$ and $\p(v,u) > 0$ for all $u \neq v$. It is also clear that $\p(v,u) = \p(u,v)$. We have left to check the triangle inequality. To verify the triangle inequality it is sufficient to show that 
\[
\fr({\norm|v+u|_n}/{1+\norm|v+u|_n}) \leq \fr({\norm|v|_n}/{1+\norm|v|_n}) +
\fr({\norm|u|_n}/{1+\norm|u|_n}).
\]
To show this, we first note that the function $f:[0,\infty) \to [0,1)$ given by $f(t) = \fr(t/1+t)$ is increasing. Thus
\begin{align*}
\fr({\norm|v+u|_n}/{1+\norm|v+u|_n}) &\leq 
\fr({\norm|v|_n + \norm|u|_n}/{1+\norm|v|_n + \norm|u|_n}) \\
&= \fr({\norm|v|_n}/{1+\norm|v|_n + \norm|u|_n}) + \fr({ \norm|u|_n}/{1+\norm|v|_n + \norm|u|_n}) \\
&\leq \fr({\norm|v|_n}/{1+\norm|v|_n}) +
\fr({\norm|u|_n}/{1+\norm|u|_n}).
\end{align*}

\begin{Prop} \label{P:metricProp}
The metric $\p$ on $V$ has the following properties:
\begin{numList}
\item $\p(v,u) = \p(v-u,0)$
\item If $v_k \to 0$ in $V$, then $\p(v_k,0) \to 0$.
\end{numList}
\end{Prop}

\begin{proof}
That $\p(v,u) = \p(v-u,0)$ for all $u,v \in V$ is obvious from the definition.

For (2), let $v_k \to 0$ in $V$. Then $\lim_{k \to \infty} \norm|v_k|_n \to 0$ for each $n$. So, for a given $\e >0$, take $N$ so that $\fr(1/2^N) < \fr(\e/2)$. Take $K$ such that for any $k > K$ we have $\norm|v_k|_n < \fr(\e/2)$ for all $1 \leq n \leq N$. Then for $k > K$ we have 
\[
\sum_{n=1}^{\infty} \fr(1/2^n) \fr({\norm|v_k|_n}/{1+\norm|v_k|_n}) 
= \sum_{n=1}^{N} \fr(1/2^n) \fr({\norm|v_k|_n}/{1+\norm|v_k|_n}) 
+ \sum_{n=N+1}^{\infty} \fr(1/2^n) \fr({\norm|v_k|_n}/{1+\norm|v_k|_n}) 
< \fr(\e/2) + \fr(1/2^N)  < \e.
\]
Therefore $\p(v_k,0) \to 0$ as $k \to \infty$.
\end{proof}

As you may have guessed, we would not take the time to talk about this metric unless it proved useful in some way. Well, it turns out that the topology induced by this metric is identical to the original topology on $V$.

\begin{Thm} \label{T:metric}
The topology on the \CNS\ $V$ induced by the metric $\p$ is equivalent to the original topology on $V$ \emph{(}i.e. the topology induced by the family of norms $\{ \norm|\cdot|_n \}_{n=1}^{\infty}$\emph{)}.
\end{Thm}

\begin{proof}
Applying \Propref{P:metricProp}, it is sufficient to consider the neighborhoods $\nBall{\delta}$ of $0$ in $V$ and the sets $\mBall$ for $\e, \delta > 0$ and $n \in \{1,2, \ldots \}$. We have to show that every $\nBall{\delta}$ contains some $\mBall$ and conversely.

Consider a neighborhood $\nBall{\delta}$ in $V$. If $v \in V$ satisfies $\p(v,0) < \e$, then $\fr(1/2^n) \fr(\norm|v|_n / 1 +\norm|v|_n) < \e$ and thus 
\[
\norm|v|_n < \fr(2^n \e / 1 - 2^n \e) = \fr( 2^n / {\fr(1/\e)} - 2^n).
\]
So, take $\e > 0$ such that 
\[
0 < \fr( 2^n / {\fr(1/\e)} - 2^n) < \delta,
\] and we have $\mBall \subset \nBall{\delta}$.

Now consider a set $\mBall$. Assume, by contradiction, there is no $n$ and $\delta > 0$ such that $\nBall{\delta} \subset \mBall$. Then for each $k$ we can find $v_k \in \nBall[k]{\fr(1/k)}$ such that $v_k$ is not in $\mBall$. This gives us a sequence $\{v_k\}_{k=1}^{\infty}$ that tends to $0$ in $V$ but not with respect to the metric $\p$. This contradicts \Propref{P:metricProp}.
\end{proof}

From this it follows that $V$ is a complete \CNS\ if and only $(V, \p)$ is a complete metric space. The following is a result which proves useful in a few theorems to come:

\begin{Lemma} \label{L:conZero}
Given a closed convex symmetric absorbing set $C$ in a complete \CNS\ $V$ we can find a neighborhood $U$ of $\,0$ contained in $C$.
\end{Lemma}

\begin{proof}
Since $C$ is absorbing we have that $V \subset \bigcup_{n=1}^{\infty} nC$. Knowing that $V$ is a complete metric space we can apply the Baire category theorem to see that the closed set $C$ is not nowhere dense. Thus the interior of $C$, $C^\circ$, is not empty. Take $v$ in $C^\circ$ and let $U$ be a symmetric open set around $0$ such that 
$v + U \subset C^\circ$ (e.g. take $U$ to be one of the $\B_k( {\fr(1/k)})$ described in \Corref{C:baseCount}). 

Because $C$ is symmetric we have that $-v - U = -v + U$ is in $C$. Since $C$ is convex it contains the convex hull of $v+U$ and $-v+U$. But this convex hull contains $U$; observe for any $w \in U$ we have that 
\[ 
w = \fr( {\left(v+w\right) + \left(-v+w\right)}/2).
\]
Thus we are done.
\end{proof}


\section{Weak Topology} \label{S:weak}

The weak topology is the simplest topology placed on the dual of a \CNS. It is defined as follows:
\begin{Def}
The \emph{weak topology} on the dual $V'$ of a \CNS\ $V$ is the coarsest vector topology on $V'$ 
such that the functional \ip<\cdot,v> is continuous for any $v \in V$. 
\end{Def}

In the following propositions, we prove some commonly used properties of the weak topology.

\begin{Prop}
The weak topology on  $V'$ has a local base of neighborhoods given by sets of the form:
\[
\N(v_1,v_2, \ldots, v_k;\e) = \weakN{\e}.
\]
\end{Prop}

\begin{proof}
In order for \ip<\cdot,v> to be continuous for all $v\in V$ we need \ip<\cdot,v> to be continuous at 0. Or equivalently, we require that $\ip<\cdot,v>^{-1}\left(-\e,\e \right) = \N(v ; \e)$ be open for each $\e \in \R$. Hence for each $v\in V$  we form the topology $\tau_v$ on $V$ given by the local base $\{ \N(v;\e) \}_{\e > 0 }$. The weak topology is the least upper bound topology for the family $\{ \tau_v \}_{v \in V}$ (see subsection   \ref{SS:famtop}). Thus, by \Thmref{T:topfam}, a local base for the weak topology is given by sets of the form
\[
\N(v_1,v_2, \ldots, v_k;\e) = \N(v_1;\e) \cap \N(v_2;\e) \cap \cdots \cap \N(v_k; \e), 
\]
where $v_1, v_2, \ldots, v_k \in V$.
\end{proof}

\begin{Prop} \label{P:wIncCont}
The inclusion map $i_n':V_n' \to V'$ is continuous when $V'$ is given the weak topology.
\end{Prop}

\begin{proof}
Consider the weak base neighborhood $\N(v_1 \dots v_k;\e)$ where $v \in V$. 
Observe that 
\[ 
i_n'^{-1} \left( \N(v_1 \dots v_k;\e) \right)= \weaknN{\e}.
\] 
Since for each $j$, $v_j \in V \subset V_n$
we have that the functional $\ip<\cdot,v_j>$ is continuous on $V_n'$. (Since $V_n$ is a Banach space, $V_n \subset V_n''$.) Thus $\weaknN{\e}$, is open in $V_n'$, being the finite intersection of open sets. 
\end{proof}

\begin{Prop}
Let $V$ be a \CHS. Then the space $V_n'$ is dense in $V'$ when $V'$ is endowed with the weak topology.
\end{Prop}

\begin{proof}
Consider $v_0' \in V'$. An arbitrary neighborhood $U$ of $v_0'$ contains a set of the form $v_0' + N$ where
$N = \N(v_1,\ldots,v_k ; \e) = \weakN{\e}$. We must find a $v_n' \in V_n'$ such that $v_n' \in v_0' + N$. 
That is $\vert \ip<v_n' - v_0', v_j> \vert < \e$ for all $ 1 \leq j \leq k$.
 
Now $v_0' \in V_l'$ for some $l$ since $V' = \bigcup_{n=1}^{\infty} V_n'$. If $l \leq n$ we are done, since $V_l' \subset V_n'$ by \Propref{P:dualInc}. If $l > n$ a little more work needs to be done, but it is still very straightforward.

 For clarity, we assume $k=2$ and $v_1, v_2$ are independent unit vectors in $V_n$. (There is no harm in assuming this. We can just shrink $\e$ suitably by dividing by the maximum of $\norm|v_1|_n$ and $\norm|v_2|_n$.) Suppose $\ip<v_0',v_1> = \l_1$ and $\ip<v_0',v_2> = \l_2$. Write $v_2$ as $v_2 = \alpha v_1 + \beta v_1^{\bot}$ where $v_1^{\bot}$ is a unit vector in the orthogonal complement of $\{ v_1 \}$ in $V_n$. Then  $\l_2= \ip<v_0',v_2> = \l_1 \alpha + \beta \ip<v_0',v_1^{\bot}>$ or equivalently $\ip<v_0',v_1^{\bot}> = \fr({\l_2 - \l_1 \alpha}/{\beta})$. Consider $w = \l_1 v_1 + \fr({\l_2 - \l_1 \alpha}/{\beta}) v_1^{\bot}$. Now $w \in V_n$. Thus $\ip<w,\cdot>_n$ is in $V_n'$, where $\ip<\cdot,\cdot>_n$ is the inner-product on $V_n$. We now observe that $\ip<w,v_1>_n = \l_1$ and $\ip<w,v_2>_n = \ip<w,\alpha v_1 + \beta v_1^{\bot}>_n = \l_1 \alpha + \l_2 - \l_1 \alpha = \l_2$. Hence $\ip<w,\cdot>_n$ agrees with $v_0'$ on $v_1$ and $v_2$.
Therefore $w \in v_0' + N$ and we have that $V_n'$ is dense in $V'$.
\end{proof}


\section{Strong Topology} \label{S:strong}

Recall the notion of bounded sets in a \CNS\ $V$ (as in subsections \ref{SS:bndV} and \ref{SS:bndR}). Using bounded sets in $V$ we can define the \emph{strong topology} on $V'$.

\begin{Def}
The \emph{strong topology} on the dual $V'$ of a \CNS\ $V$ is defined to be the topology with a local base given by sets of the form
\begin{center}
\(
\N(D;\e) = \strongN{\e},
\)
\end{center}
where $D$ is any bounded subset of $V$ and $\e > 0$.
\end{Def}

Taking $D$ to be a finite set such as $\{ v_1, v_2, \ldots, v_k \}$, it is clear that the strong topology is finer than the weak topology.

\begin{Prop} \label{P:sIncCont}
The inclusion map $i_n':V_n' \to V'$ is continuous when $V'$ is given the strong topology.
\end{Prop}

\begin{proof}
Consider the neighborhood $\rN(D;\e) = \strongN{\e}$ where $D$ is a bounded set in $V$ and $\e > 0$. 
Now 
\[
i_n'^{-1}(\N(D;\e)) = \bstrongnN{\e}.
\] 
Let $\sup_{v \in D} \norm|v|_n = M$. Take $v_0'$ in $i_n'^{-1}(\rN(D;\e))$ and let $c_0 = \sup_{v \in D} \vert \ip<v_0',v> \vert < \e$. 
Consider the open set $B(v_0',\fr({\e-c_0}/M+1 )) = \{ v' \in V_n' ; \medspace \norm|v'-v_0'|_{-n} < \fr({\e - c_o}/M+1) \}$.
We assert that $B(v_0',\fr({\e-c_0}/M+1)) \subset i_n'^{-1}(\rN(D;\e))$. 

Take $v' \in B(v_0',\fr({\e-c_0}/M+1))$. Then $\norm|v' - v_0'|_{-n} < \fr({\e - c_0}/M+1)$. 
This gives us the following 
\[
\sup_{v \in D} \big|\ip<v' - v_0', {\tfr(v/{\norm|v|_n})}> \big| < \fr({\e - c_0}/M+1).
\]
Thus $\sup_{v \in D} | \ip<v' - v_0',v>| < \e - c_0$, since $\norm|v|_n \leq M$ when $v \in D$. 
From this we see that $\sup_{v \in D}  |\ip<v',v>| < \e$. Therefore $v' \in i_n'^{-1}(\N(D;\e))$.
\end{proof}

\subsection{Strongly bounded sets of $V'$}
When $V'$ is endowed with the strong topology, a bounded set $B \subset V'$ is called \emph{strongly bounded}. (Likewise when $V'$ has the weak topology, $B$ is said to be \emph{weakly bounded}). Strongly bounded sets have many nice properties, which we will prove in this section. First let us begin with the following definition:

\begin{Def}
A set $B \subset V'$ is said to be \emph{bounded on the set} $A \subset V$ if
\[
\sup_{v' \in B, v \in A} \left| \ip<v',v> \right| < \infty.
\]
\end{Def}

\begin{Lemma} \label{L:bndIff}
A set $B \subset V'$ is strongly bounded if and only if it is bounded on each bounded set $D \subset V$.
\end{Lemma}

\begin{proof} ($\Rightarrow$)
Let $B \subset V'$ be strongly bounded and let $D$ be a bounded set of $V$. Consider the neighborhood of $V'$ given by 
\[
\N(D;1) = \dstrongN{1}.
\]
Since $B$ is bounded there exists an $\l > 0$ such that $B \subset \l \rN(D;1)$ 
or equivalently $\tfr(1/{\l}) B \subset\rN(D;1)$. Then for $v' \in B$ we have $\fr(v'/{\l}) \in \rN(D;1)$.
Thus $|\ip<v,v'>| \leq \l$ for any $v \in D$. Therefore $B$ is bounded on the set $D$.

$(\Leftarrow)$
Suppose $B$ is bounded on each bounded set $D \subset V$. Consider a neighborhood $\N(D;\e)$ of $0$ in $V'$.
By hypothesis, $\sup_{v' \in B, v \in D} |\ip<v',v>| = M < \infty$. So for any $v' \in B$ we have that
$|\ip<\fr(\e v'/M+1), v>| < \e$ when $v \in D$. Thus $\fr(\e/M+1) B \subset \rN(D;\e)$ or equivalently $B \subset \fr(M+1/{\e}) \rN(D;\e)$. Hence $B$ is bounded.
\end{proof}

\begin{Lemma} \label{L:bndBk}
A set $B \subset V'$ is strongly bounded if and only if there exists $k$ such that $B$ is bounded on $B_k({\fr(1/k)})$.
\end{Lemma}

\begin{proof} ($\Rightarrow$)
As per \Corref{C:baseCount}, consider the local base sets $\rB_1(1) \supset \rB_2({\tfr(1/2)}) \supset \cdots$ of $V$. By contradiction suppose that $B$ is not bounded on $B_k({ \fr(1/k)} )$ for any $k$. Then for every $k$ there exist a $v_k \in B_k({\fr(1/k)})$ and a $v_k' \in B$ such that 
$| \ip<v_k', v_k>| > k $. The sequence $\{v_k \}$ goes to 0, thus it must be bounded. So by \Lemmaref{L:bndIff} there must exist a positive number $M$ such that 
$|\ip<v',v_k>| \leq M$ for all $v' \in B$ and all $k \in \{1,2,\ldots \}$. This contradicts the way $v_k'$ and $v_k$ were chosen.

($\Leftarrow$)
Conversely, let $B \subset V'$ be bounded on some $B_k({ \fr(1/k)} ) \subset V$. Take a bounded set $D$ in $V$. Then $D \subset \l \B_k({\fr(1/k)})$ for some $\l > 0$. Thus $B$ is bounded on $D$, since $B$ is bounded on $\l \B_k({\fr(1/k)})$. Thus by \Lemmaref{L:bndIff}, $B$ is bounded.
\end{proof}

\begin{Thm} \label{T:bndIff}
A set $B \subset V'$ is strongly bounded if and only if $B \subset V_k'$ for some $k$ and $B$ is bounded in the norm $\norm|\cdot|_{-k}$ of $V_k'$.
\end{Thm}

\begin{proof}
($\Leftarrow$)
Let $B \subset V_k'$ be bounded in the norm $\norm|\cdot|_{-k}$ by some  $M>0$ (i.e. $\sup_{v' \in B} \norm|v'|_{-k} < M$). Consider the set $\rB_{k}(1) = \{ v \in V ; \norm|v|_k < 1 \}$.
Then for $v' \in B$ and $v \in \rB_{k}(1)$ we have that $|\ip<v',v>| \leq M$. Thus $B$ is bounded on $\B_k(1)$ and hence on $\rB_k( {\fr(1/k)} ) $. Therefore $B$ is strongly bounded by \Lemmaref{L:bndBk}.

($\Rightarrow$) Conversely suppose $B$ is a strongly bounded set in $V'$. Then by \Lemmaref{L:bndBk} there is a $k$ such that $B$ is bounded on the set $\rB_k( {\fr(1/k)} ) = \{ v \in V ; \norm|v|_k < \tfr(1/k) \}$. That is there is an $M < \infty$ such that 
$ |\ip<v',v>| \leq M$  for all  $v' \in B$ and all  $v \in \B_k( {\fr(1/k)})$.

Let $N_k \subset V_k$ be given by $N_k = \{ v \in V_k ; \norm|v|_k < \fr(1/k) \}$.
Since $V$ is dense in $V_k$ we have that 
\[
\sup_{v' \in B, v \in N_k} |\ip<v',v>| \leq M.
\] 
From the above we see for any $v' \in B$ and unit vector $v \in V_k$ we have that 
$|\ip<v', \fr(v/ {\left( k+1 \right)} ) >| < M$. Hence $|v'|_{-k} \leq \left( k+1 \right) M$. 
Thus for any $v' \in B$ we have that $v' \in V_k'$ and $\norm|v'|_{-k} \leq \left( k+1 \right) M$.
\end{proof}

\subsection{Reflexivity}
 Just as we can discuss the dual $V'$ of $V$, we can also talk about the dual of $V'$. Of course, this depends on the topology we put on $V'$. As we will see it turns out that $V'' = V$ as sets if $V'$ is given the weak or strong topology (and $V$ is a \CHS). We can also put a topology on $V''$. We construct this topology from the strongly bounded sets in $V'$. For each set $B$ in $V'$ that is strongly bounded and each $\e > 0$ form the neighborhood 
\[
\N(B;\e) = \ddualStrongN{\e}.
\]
Take the collection of all sets $\rN(B;\e)$ as our local base in $V''$. We call this topology the \emph{strong topology on} $V''$. Given this topology we will also see that $V''$ is homeomorphic to $V$.

\begin{Prop} \label{P:refl}
Let $V$ be a \CHS. Then $V = V''$ when $V'$ is given the weak or strong topology.
\end{Prop}

\begin{proof}
Consider $v \in V$ and the corresponding linear functional $\hat{v}$ on $V'$ given by
\[
\ip<\hat{v}, v'> = \ip<v',v>, \text{\qquad where } v' \in V'.
\] 
Observe that $\ip<\hat{v}, \cdot>$ is continuous since 
$\ip<\hat{v}, \cdot>^{-1} \left( -\e, \e \right) = \{ v' \in V' ; |\ip<v',v>| < \e \}$ 
which is open in the weak (and hence the strong) topology on $V'$.

Also note that if $\hat{u} = \hat{v}$, then $\ip<v',v> = \ip<v', u>$ for all $v' \in V'$. Thus $v = u$. Therefore the correspondence $v \to \hat{v}$ is injective.

We now show that the correspondence $v \to \hat{v}$ is surjective. Take $v'' \in V''$. 
Then $v''$ is continuous on $V'$. Since, by \Propref{P:dualInc}, $V' = \bigcup_{n=1}^{\infty} V_n'$  we have that $v'' \in V_n''$ for all $n$. But $V_n = V_n''$ since $V_n$ is a Hilbert space. Thus $v''$ can be considered as an element of $V_n$ for all $n$. Since $V$ is a \CHS\ we have that $\cap_{n=1}^{\infty} V_n = V$ by \Propref{P:complete}. Thus $v'' \in V$ and we have that 
$v \to \hat{v}$ is surjective.
\end{proof}

\begin{Thm}
If $V$ is a \CHS, then $V''$ is homeomorphic to $V$ when $V''$ is given the strong topology.
\end{Thm}

\begin{proof}
From \Propref{P:refl} we already see that $V = V''$. We now need to see that the correspondences $\hat{v} \to v$ and $v \to \hat{v}$ are continuous.

First we consider the continuity of $v \to \hat{v}$. Let $\N(B;\e)$ be a neighborhood of $0$ in $V''$.
So we have that $B$ is a strongly bounded set in $V'$. By \Thmref{T:bndIff} we know that
$B \subset V_k'$ for some $k$ and is bounded in the norm $\norm|\cdot|_{-k}$. 
Let us call $\sup_{v' \in B} \norm|v'|_{-k} = M < \infty$. Consider the neighborhood $\B_k({\fr(\e/M)}) \subset V$ given by 
$\B_k({\fr(\e/M)}) = \{ v \in V \medspace ; \medspace \norm|v|_k < \fr(\e/M) \}$. Take a $v \in \B_k({\fr(\e/M)})$. We need to see that $\hat{v} \in \N(B;\e)$. For any $v' \in B$ we have that 
\[
\aip<\hat{v},v'> = \aip<v',v> \leq \norm|v'|_{-k} \norm|v|_k < M \fr(\e/M) = \e.
\]
So $\hat{v} \in \N(B;\e)$. Thus $v \to \hat{v}$ is continuous.

Now consider $\hat{v} \to v$. Let $0 < \e < 1$ and take $\B_k(\e) = \{ v \in V ; \norm|v|_k < \e \}$, a member of the local base for $V$ (see subsection \ref{SS:openV}). 
Let $B \subset V'$ be given by 
\[
B = \{ v' \in V' ; |\ip<v', v>| \leq 1, \text{ for all } v \in V_k \text{ with } \norm|v|_k < \e \}.
\]
Note that $B$ is strongly bounded by \Thmref{T:bndIff}. So we can form the local base element 
$\N(B;\e)$ of $V''$ given by
\[
\N(B;\e) = \dualStrongN{\e}.
\]
Take a $\hat{v} \in \N(B;\e)$. Note that $\ip<\tfr(v/\norm|v|_k), \cdot>_k \in B$ since $\ip<\tfr(v/\norm|v|_k), u>_k \leq \norm|u|_k$ for $u \in V_k$. Since  $\hat{v} \in \N(B;\e)$ and $\ip<\tfr(v/\norm|v|_k), \cdot>_k \in B$, we must have
$|\ip<\tfr(v/\norm|v|_k), v>_k| = \norm|v|_k < \e$.
Therefore $v \in \rB_k(\e)$.  This proves the continuity of the map $\hat{v} \to v$.
\end{proof}

\subsection{Completeness in $V'$} 
Suppose $V'$ is given the strong topology.  The convergence of a sequence of functionals $\{v_k'\}_{k=1}^{\infty}$ in $V'$ to an element $v_0' \in V'$ is called \emph{strong convergence} and $\{v_k'\}_{k=1}^{\infty}$ is said to \emph{converge strongly} to $v_0'$. Obviously $\{v_k'\}_{k=1}^{\infty}$ converging strongly to $v_0'$ is equivalent to $\{v_k' -v_0'\}_{k=1}^{\infty}$ converging strongly to $0$. Thus a sequence $\{v_k'\}_{k=1}^{\infty}$ converges strongly to $v_0'$ if and only if for any bounded set $D \subset V$ and any number $\e > 0$ there exists a $K > 0$ such that $v_k' - v_0' \in \N(D;\e) = \strongN{\e}$ for all $k \geq K$. Hence a sequence $\{v_k'\}_{k=1}^{\infty}$ converges strongly to $v_0'$ if and only if 
$\{\ip<v_k',\cdot> \}_{k=1}^{\infty}$ converges uniformly to $\ip<v_0',\cdot>$ on each bounded set $D \subset V$. We say that a sequence $\{v_k'\}_{k=1}^{\infty}$ is \emph{strongly Cauchy} (or \emph{strongly fundamental}) if the sequence of numbers $\{\ip<v_k', v>\}_{k=1}^{\infty}$ converges for each element $v \in V$ and the convergence is uniform on each bounded set $D \subset V$.

\begin{Thm}
Let $V$ be a \CNS. The dual $V'$ of $V$ is complete under the strong topology.
\end{Thm}

\begin{proof}
Let $\{v_k'\}_{k=1}^{\infty}$ be a strongly Cauchy sequence in $V'$. Then for $v \in V$ we have that the sequence of numbers  $\{\ip<v_k',v>\}_{k=1}^{\infty}$ converges. We conveniently denote this limit by $\ip<v',v>$. For each $v \in V$ we have
\[
\ip<v',v> = \lim_{k \to \infty} \ip<v_k',v>.
\]

This functional $\ip<v',\cdot>$ is clearly linear on $V$. We have to check that it is continuous. 
For this it is sufficient to see that $\ip<v',\cdot>$ is bounded on bounded sets (see \Propref{P:borno}).
Let $D$ be a bounded set in $V$. Observe that the functions $\{\ip<v_k',\cdot>\}_{k=1}^{\infty}$ are bounded on $D$. Moreover they converge uniformly to $\ip<v',\cdot>$ on $D$. Hence there is a $K > 0$ such that $\aip<v' -v_K',v> < 1$ for all $v$ in $D$. Thus we have that 
\[
\sup_{v \in D} \aip<v',v> \leq \sup_{v \in D} \aip<v_K',v> + 1 < \infty.
\]
Therefore $\ip<v',\cdot>$ is bounded on bounded sets and hence continuous. So $v' \in V'$ and $V'$ is complete with respect to the strong topology.
\end{proof}

\subsection{Comparing the Weak and Strong topology}
When a \CNS\ $V$ is complete, many properties of the strong and weak topologies coincide. We will see that weakly and strongly bounded sets are one in the same. Also under suitable conditions, weak and strong convergence coincide.

\begin{Thm} \label{T:wsbnd}
Let $V$ be a complete \CNS\ with dual $V'$. Every weakly bounded set in $V'$ is strongly bounded.
\end{Thm}

\begin{proof}
By \Lemmaref{L:bndBk} and \Corref{C:baseCount} it is sufficient to show that a weakly bounded set $B$ is bounded on some neighborhood of zero in $V$.

Let us define a set $C \subset V$ as follows:
\[
C = \{ v \in V ; \aip<v',v> \leq 1 \text{ for all } v' \in B \}
  = \bigcap_{v' \in B} \{ v \in V ; \aip<v',v> \leq 1 \}.
\]
Observe that $C$ is closed, being the intersection of closed sets, $C$ is convex, being the intersection of convex sets, and $C$ is symmetric, being the intersection of symmetric sets. Finally note that $C$ is absorbent: Take $v \in V$. Since $B$ is weakly bounded we must have $B \subset \l \N(v;1)$ where $\N(v;1) = \weakS{1}$ for some $\l > 0$. Thus
$\aip<v',v> \leq \l$ for all $v' \in B$. Hence $\fr(v/{\l})  \in C$ or equivalently $v \in \l C$.

So we can apply \Lemmaref{L:conZero} to see that there is a neighborhood $U$ of $0$ in $V$ such that $U \subset C$.  Therefore the elements of $B$ are uniformly bounded on $U$ (by $1$). Thus $B$ is bounded on $U$ and hence strongly bounded.
\end{proof}

\begin{Cor} \label{C:pwBnd}
Let $V$ be a complete \CNS\ with dual $V'$. If a sequence $\{ v_k' \}_{k=1}^{\infty}$ in $V'$ converges pointwise $($on each $v \in V)$, then $\{ v_k' \}_{k=1}^{\infty}$ is strongly bounded.
\end{Cor}
\begin{proof}
Since $\{ v_k' \}_{k=1}^{\infty}$ converges pointwise, it is weakly bounded.
\end{proof}

\begin{Cor}
Let $V$ be a complete \CNS\ with dual $V'$. Then $V'$ is complete with respect to the weak topology.
\end{Cor}

\begin{proof}
Take a Cauchy sequence $\{ v_k' \}_{k=1}^{\infty} \subset V'$. Then by \Corref{C:pwBnd}, we have that $\{ v_k' \}_{k=1}^{\infty}$ is strongly bounded. Thus by \Lemmaref{L:bndBk}, 
$\{ v_k' \}_{k=1}^{\infty}$ is bounded on some neighborhood $U$ of $0$ in $V$. That is, there exists an $M > 0$ such that $\aip<v_k',v> \leq M$ for all $v\in U$ and all $k \in \{ 1,2, \ldots 
\}$.

Define $v'$ by $\ip<v',v> = \lim_{k \to \infty} \ip<v_k',v>$. Obviously, $v'$ is linear. Observe that for all $v \in U$ we have
\[
\aip<v',v> = \lim_{k \to \infty} \aip<v_k',v> \leq M.
\]
So $v'$ is bounded on $U$ and hence continuous (by \Corref{C:borno}).
\end{proof}

Of particular interest are \CNS s such with the property that bounded sets are limit point compact (see \defref{D:lpc}). These spaces have many wonderful properties, that do not hold in general for infinite-dimensional normed spaces. We make the following definition (the terminology comes from Gel'fand \cite{iG68}):

\begin{Def} 
A complete \CNS\ $V$  in which all bounded sets are limit point compact is called \emph{perfect}.
\end{Def}

\begin{Remark}
Since \Thmref{T:metric} gives us that $V$ is metrizable,  limit point compact can be replaced with compact or sequentially compact in the above definition. Therefore if $V$ is perfect, the strong topology on $V'$ is nothing more than the well known \emph{compact--open topology} \cite{jM00}. 
\end{Remark}

\begin{Thm}
Let $V$ be a perfect space with dual $V'$. Then
a sequence $\{ v_k' \}_{k=1}^{\infty}$ in $V'$ converges strongly if and only if it converges weakly. \emph{(}i.e. weak and strong converge coincide on the dual space $V')$
\end{Thm}
\begin{proof}
Obviously strong convergence implies weak convergence. So take a sequence $\{ v_k' \}_{k=1}^{\infty}$ in $V'$ which converges weakly to $v' \in V'$. Without loss of generality we can take $v' = 0$ (replace $v_k'$ with $v_k' - v'$). The sequence 
$\{ v_k' \}_{k=1}^{\infty}$ is weakly bounded, being weakly convergent. Thus by \Thmref{T:wsbnd} we have that $\{ v_k' \}_{k=1}^{\infty}$ is strongly bounded.

To show $\{ v_k' \}_{k=1}^{\infty}$ converges strongly we must show $\ip<v_k',\cdot>$ goes to $0$ uniformly on each bounded set $D \subset V$. Suppose, by contradiction, that there exists a bounded set $D$ in $V$ where $\ip<v_k',\cdot>$ does not go to $0$ uniformly. So for some $\e > 0$, there is a $k_1 \geq 1$ such that $\aip<v_{k_1}',v> \geq \e$ for some $v \in D$. Name this $v$ as $v_{k_1} $. Likewise there an $k_2 > k_1$ and $v_{k_2} \in V$ such that $\aip<v_{k_2}',v_{k_2}> \geq \e$. Continuing in this manner we form a sequence $\{ v_{k_j} \}_{j=1}^{\infty}$ in $V$. This sequence is bounded, being taken from the bounded set $D$.

Knowing that $V$ is a perfect space we have a subsequence of $\{ v_{k_j} \}_{j=1}^{\infty}$ that converges to some $v$ in $V$. Renumbering if necessary we will just take this subsequence to be $\{ v_{k_j} \}_{j=1}^{\infty}$.  Since $\{ v_{k_j} \}_{j=1}^{\infty}$ goes to $v$ in $V$ then the sequence given by $w_{k_j} = v_{k_j} - v$ goes to $0$ in $V$. 

Now for any strongly bounded set $B \subset V'$, \Thmref{T:bndIff} guarantees that $\ip<v',w_{k_j}>$ goes to $0$ uniformly for all $v' \in B$. So take $B$ to be the set
$\{ v_{k_j}' \}_{j=1}^{\infty}$. Then $\ip<v_{k_j}',w_{k_j}>$ goes to $0$ and by weak convergence we have that $\ip<v_{k_j},v>$  goes to $0$. Thus
\[
\lim_{j \to \infty} \ip<v_{k_j}', v_{k_j}> = 
\lim_{j \to \infty} \ip<v_{k_j}',w_{k_j}> + \ip<v_{k_j}',v> = 0.
\]
This contradicts the construction of the $v_{k_j}$ and $v_{k_j}'$.
\end{proof}


\section{Inductive Limit Topology} \label{S:induct}

Given a sequence of normed spaces $\{ (W_n, \norm|\cdot |_n ); n \geq 1 \}$ with $W_n$ continuously imbedded in $W_{n+1}$ for all $n$, we form the space $W = \bigcup_{n =1}^{\infty} W_n$ and endow $W$ with the finest locally convex vector topology such that for each $n$ the inclusion map $i_n:W_n \to W$ is continuous. This topology is call the \emph{inductive limit topology} on $W$ and $W$ is said to be the \emph{inductive limit} of the sequence $\{ (W_n, \norm|\cdot |_n ); n \geq 1 \}$

\subsection{Local Base} As always, when discussing a vector topology, we should try to discover what a useful local base for the topology would be.

\begin{Thm} \label{T:baseIT}
Suppose $W$ is the inductive limit of the normed spaces $\{ (W_n, \norm|\cdot |_n ); n \geq 1 \}$. A local base for $W$ is given by the set \sB\ of all balanced convex subsets $U$ of $W$ such that $i_n^{-1}\left( U \right)$ is a neighborhood of $0$ in $W_n$ for all $n$.
\end{Thm}

\begin{proof}
We first apply \Lemmaref{L:convexBase} to see the set \sB\ is in fact a local base for $W$. Take $U,V \in \sB$, then clearly $U \cap V \in \sB$. Now if $U \in \sB$, then it is easy to see that $\alpha U \in \sB$ for $\alpha \neq 0$. Finally we show $U \in \sB$ is absorbing. Note that $i_n^{-1}(U) $ is absorbing in $W_n$ (since $W_n$ is a normed space and $i_n^{-1}(U) $ is open in $W_n$). Thus $U$ absorbs all the point of $W_n = i_n(W_n) \subset W$. Since $ W = \bigcup_{n=1}^{\infty} W_n$, $U$ absorbs $W$. Thus by \Lemmaref{L:convexBase} we see that $\sB$ is a base of neighborhoods for a locally convex vector topology on $W$.

It is fairly straightforward to see that \sB\ gives us the finest locally convex vector topology making all the $i_n:W_n \to W$ continuous: Let $\tau$ be a locally convex vector topology on $W$ making all the $i_n$ continuous. Take a convex neighborhood (of 0) $U$ in $\tau$. By \Lemmaref{L:balance} we can assume $U$ is balanced. Since each $i_n$ is continuous, we have $i_n^{-1}(U) $ is a neighborhood in $W_n$. Thus $U \in \sB$.
\end{proof}

\begin{Cor} \label{C:baseIT}
Suppose $W$ is the inductive limit of the normed spaces $\{ (W_n, \norm|\cdot |_n ); n \geq 1 \}$. A local base for $W$ is given by the balanced convex hulls of sets of the form $\bigcup_{n=1}^{\infty}i_n({\rB_n(\e_n)})$ \emph{(}where 
$\rB_n(\e_n) = \Wnball{\e_n}$\emph{)}.
\end{Cor}

\begin{proof} 
Let $U$ be the balanced convex hull of the set $\bigcup_{n=1}^{\infty}i_n({\rB_n(\e_n)})$ in $W$. Then $\B_n(\e_n) \subset i_n^{-1}(U)$. So $i_n^{-1}(U)$ is a neighborhood of $0$ in $W_n$. By \Thmref{T:baseIT} such a $U$ is a neighborhood in $W$.

Now if $U$ is any balanced convex neighborhood of $0$ in $W$, then $i_n^{-1}(U)$ contains a neighborhood $\rB_n(\e_n)$. Hence $i_n(\rB_n(\e_n)) \subset U$. Since $U$ is convex and balanced, the balanced convex hull of $\bigcup_{n=1}^{\infty}i_n({\rB_n(\e_n)})$ is contained in $U$. Thus the sets described form a local base for $W$.
\end{proof}

\subsection{Inductive Limit Topology on $V'$} \label{SS:iltOnDaul}
Let $V$ be a \CNS . Then $V'$, the dual of $V$, can be regarded as the inductive limit of the sequence of normed spaces $\{(V_n', \norm|\cdot|_{-n}); n \geq 1 \}$. Thus $V'$ can be given the inductive limit topology. In light of \Propref{P:wIncCont} and \Propref{P:sIncCont} we see that the inductive limit topology on $V'$ is finer than the strong and weak topology on $V'$. We also have the following useful result about convergence on $V'$ in the inductive topology:

\begin{Thm}
Let $V$ be a \CNS. Endow $V'$ with the inductive limit topology. A sequence $\{v_k'\}_{k=1}^{\infty}$ converges to $v'$ in $V'$ if and only if there exists some $n$ such that $v_k \in V_n'$ for all $k$ and $\lim_{k\to \infty} \norm|v_k'-v|_{-n} = 0$ $($i.e. $v_k'$ converges to $v'$ in $V_n')$.
\end{Thm}
\begin{proof}
($\Leftarrow$) Using \Corref{C:baseIT}, this direction is obvious.

($\Rightarrow$) Let $\{v_k'\}_{k=1}^{\infty}$ be a sequence in $V'$ that converges to $v' \in V'$. Replacing $v_k'$ with $v_k' -v'$, if necessary, we assume that $v' = 0$. Since $\{v_k'\}_{k=1}^{\infty}$ converges to $0$ in the inductive limit topology, by the above discussion, it converges to $0$ in the strong topology on $V$. Hence $\{v_k'\}_{k=1}^{\infty}$ is strongly bounded. Thus by \Thmref{T:bndIff} we have that there is an $n$ such that $\{v_k'\}_{k=1}^{\infty} \subset V_n'$.

Now we must show that $\norm|v_k'|_{-n}$ goes to $0$ as $k$ tends to infinity. That is for a given $\e > 0$ we need to find a $K > 0$ such that for all $k \geq K$ we  have $\norm|v_k'|_{-n} < \e$. Consider the base neighborhood $U$ of $V'$ given by the balanced convex hull of $\bigcup_{l=1}^{\infty} B_l$ where for $l = n$ we take
\[
B_n = \{v' \in V_n' \, ;\, \norm|v'|_{-n} < \e \}.
\]
For $l < n$, $B_l = \{v' \in V_l' \, ;\, \norm|v'|_{-l} < \e_l \}$ where $\e_l > 0$ is chosen so that $B_l$ is contained in $i_{l,n}'^{-1}(B_n)$. 
(Such an $\e_l > 0$ exist by the continuity of the inclusion map $i_{l,n}':V_l' \to V_n'$.) 

And for $l > n$ we first note that the restricted inclusion $\tilde{i}_{n,l}:V_n' \to i_{n,l}(V_n')=V_n' \subset V_l'$ is a homeomorphism (since $V_n'$ is continuously imbedded into $V_{n+1}'$ for each $n$). This gives us $i_{n,l}'(B_n) \cap V_n' = W\cap V_n'$ where $W$ is open in $V_l'$. Thus take $B_l = \{ v' \in V_l' \, ;\, \norm|v'|_{-l} < \e_l \}$ where $\e_l > 0$ is chosen so that $B_l \subset W$. 

Now since $U$ is open, there is a $K$ such that for all $k \geq K$ we have that $v_k' \in U$. We will show that $v_k' \in B_n$ for $k \geq K$. Let $k \geq K$ and consider the element $v_k'$. Since $v_k' \in U$ we can write $v_k'$ as $v_k' = \sum_{j = 1}^{m} \l_j y_j$ where $\sum_{j=1}^{m} |\l_j| \leq 1$ and $y_j \in B_j$. Observe that each $y_j$ with $\l_j \neq 0$ is in $V_n'$. (If there is an $y_j$ not in $V_n'$ with $\l_j \neq 0$, then $v_k'$ could not be $V_n'$.) Thus we have
\begin{equation} \label{E:lt}
\norm|v_k'|_{-n} \leq \sum_{j=1}^{m} |\l| \norm|y_j|_{-n}.
\end{equation}

Observe that for $j \leq n$, $y_j \in B_j \subset i_{j,n}'^{-1}(B_n)$. So $\norm|y_j|_{-n} < \e$. Also for $j > n$ we have that $y_j \in B_j$. Since $y_j \in V_n'$ we get that $y_j \in B_j \cap V_n' \subset i_{n,j}'(B_n)\cap V_n'$. So $\norm|y_j|_{-n} < \e$. 

Therefore in \eqref{E:lt} we have that 
\[
\norm|v_k'|_{-n} \leq \sum_{j=1}^{\infty} |\l| \norm|y_j|_{-n} < \sum_{j=1}^{\infty} |\l_j| \e \leq \e.
\]
Thus $v_k'$ is in $B_n$ for all $k \geq K$ and we are done.

\end{proof}


\section{Comparing the Three Topologies}
In this section we compare the three topologies on the dual $V'$ of a \CNS\ $V$. In order to do this efficiently we first introduce a fourth topology on $V'$. It is the \emph{Mackey topology} on $V'$.

\subsection{Mackey Topology} In order to talk about the Mackey topology we need the following notion:

\begin{Def}
Let \D\ be a set of bounded subsets of a topological vector space $E$ with dual $E'$. The topology of \emph{uniform convergence on the sets of} $\D$ is the topology with subbasis neighborhoods of $0$ given by 
\[
\N(D;\e) = \Bigl\{v' \in E' \, ; \, \sup_{v \in D} \aip<v',v> < \e \Bigr\},
\]
where $D \in \D$ and $\e >0$.
This is also referred to as the topology of \emph{\D--convergence} on $E'$.
\end{Def}

From the definition we see that a local base neighborhood for the topology of  \D--convergence on a vector space $E$ with dual $E'$ looks like
\[
\N(D_1;\e_1) \cap \N(D_2; \e_2) \cap \cdots \N(D_k; \e_k),
\]
where $D_j \in \D$ and $\e_j > 0$ for all $1 \leq j \leq k$.
We now state the following theorem without proof:

\begin{Thm}[Mackey-Arens] \label{T:MA}
Suppose that under a locally convex vector topology $\tau$, $E$ is a Hausdorff space. Then $E$ has dual $E'$ under $\tau$ if and only if $\tau$ is a topology of uniform convergence on a set of balanced convex weakly--compact subsets of $E'$.
\end{Thm}

For a proof of this results see \cite{aR64}, \cite{hS66}, or \cite{gK69}. Using this theorem we can define the Mackey topology as follows:

\begin{Def}
Let $E$ be a topological vector space with dual $E'$. The \emph{Mackey topology} on $E$ is the topology on uniform convergence on all balanced convex weakly--compact subsets of $E'$.
\end{Def}

\begin{Remark} \label{R:mackeyDual}
From this discussion we see that the Mackey topology on $V'$ has a local base given by 
\[
\N(C;\e) = \Bigl\{ v' \in V' \, ; \, \sup_{v \in C} \aip<v',v> < \e \Bigr\},
\]
where $\e > 0$ and $C$ is a balanced convex weakly--compact set in $V$.
\end{Remark}

\begin{Remark} \label{R:weakCompact}
Although we have not defined the term weakly--compact, it is nothing to fret about. Just as we have defined the weak topology on $V'$, we can define an analogous topology on $V$. This topology has as its local base sets of the form
\[
\N(v_1',v_2',\ldots,v_k';\e) = \weakNV{\e}.
\]
When a set in $V$ is said to be weakly--compact, it simply means that the set is compact with respect to the weak topology on $V$.
\end{Remark}

\subsection{The Topologies on $V'$} Let us make the following notational convention throughout this section:
\begin{Notation}
Let $V$ be a \CNS\ with dual $V'$. The weak topology, strong topology, inductive limit topology, and Mackey topology on $V'$ with be denoted by $\gt_w$, $\gt_s$, $\gt_i$, and $\gt_m$, respectively.
\end{Notation}

\begin{Prop} \label{P:wCompWBnd}
Let $V$ be a \CNS. Suppose $C \subset V$ is weakly--compact, then $C$ is weakly bounded.
\end{Prop}

\begin{proof}
Let $v' \in V'$ and $\e > 0$ be given. We have to show there exists a $k$ such that $C \subset k \N(v';\e)$ (see Definition~\ref{D:weakBnd}). Cover $C$ by the sets $\{k \N(v';\e) \}_{k=1}^{\infty}$. Since $C$ is weakly--compact, $C \subset k \N(v';\e)$ for some $k$.
\end{proof}

\begin{Cor} \label{C:sFinerM}
Let $V$ be a \CNS\ with dual $V'$. Then the strong topology $\gt_s$ is finer than the Mackey topology $\gt_m$ on $V'$.
\end{Cor}

\begin{proof}
The topology $\gt_s$ is by definition the topology of uniform convergence on all bounded sets in $V$. But by \Thmref{T:bndEqV} every bounded set in $V$ is weakly bounded. And by \Propref{P:wCompWBnd}, we have that every weakly--compact set is weakly bounded. Thus $\gt_m \subset \gt_s$.
\end{proof}

\begin{Lemma} \label{L:dualHaus}
Let $V$ be a \CNS\ with dual $V'$. Then $V'$ is Hausdorff in the weak topology $t_w$, and hence in the strong, Mackey, and inductive limit topologies.
\end{Lemma}

\begin{proof}
Take $u' \in V'$. We must find a neighborhood of $0$ in $\gt_w$ that does not contain $u'$. Take $v \in V$ such that  $\aip<u',v> \neq 0$. Let $\aip<u',v> = \l \neq 0$. Consider the set $\rN(v; {\fr(\l/2) }) = \weakS{\fr(\l/2)}$. This set cannot contain $u'$. Thus $V'$ is Hausdorff in the weak topology (and hence in the finer strong, Mackey, and inductive topologies).
\end{proof}

\begin{Lemma} \label{L:dual}
Let $V$ be a \CHS\ with dual $V'$. Then the dual of $V'$ is $V$ when $V'$ is given the weak, strong, Mackey, or inductive limit topology.
\end{Lemma}

\begin{proof}
Consider $v \in V$ and the corresponding linear functional $\hat{v}$ on $V'$ given by
\[
\ip<\hat{v}, v'> = \ip<v',v> \text{\qquad where } v' \in V'
\] 
Observe that $\ip<\hat{v}, \cdot>$ is continuous since 
$\ip<\hat{v}, \cdot>^{-1}\! \left( -\e, \e \right) = \{ v' \in V' ; |\ip<v',v>| < \e \}$ 
which is open in the weak topology (and hence the strong, Mackey, and inductive limit topologies) on $V'$.

Also note that if $\hat{u} = \hat{v}$, then $\ip<v',v> = \ip<v', u>$ for all $v' \in V'$. Thus $v = u$. Therefore the correspondence $v \to \hat{v}$ is injective.

We now show that the correspondence $v \to \hat{v}$ is surjective. Take $v'' \in V''$, the dual of $V'$. 
Then $v''$ is continuous on $V'$. Since $V' = \bigcup_{n=1}^{\infty} V_n'$ by \Propref{P:dualInc}, we have that $v'' \in V_n''$ for all $n$. But $V_n = V_n''$ since $V_n$ is a Hilbert space. Thus $v''$ can be considered as an element of $V_n$ for all $n$. Since $V$ is a complete we have that $\bigcap_{n=1}^{\infty} V_n = V$ by \Propref{P:complete}. Thus $v'' \in V$ and we have that 
$v \to \hat{v}$ is surjective.
\end{proof}

\begin{Thm}
Let $V$ be a  \CHS\ with dual $V'$. Then the inductive, strong, and Mackey topologies on $V'$ are equivalent \emph{(}i.e. $\gt_s = \gt_i = \gt_m $\emph{)}.
\end{Thm}

\begin{proof}
By \Lemmaref{L:dualHaus} and \Lemmaref{L:dual} we have that $V'$ is Hausdorff and has dual $V$ under the topologies $t_s$, $t_i$, and $t_m$. Thus we can apply \Thmref{T:MA} to $V'$. (In the theorem we are taking $V'$ as $E$ and $V$ as $E'$.) This gives us that $t_i$ and $t_s$ are topologies of uniform convergence on a set of balanced convex weakly--compact subsets of $V$. However, by definition, the Mackey topology, $\gt_m$, is the finest such topology. Thus $\gt_s \subset \gt_m$ and $\gt_i \subset \gt_m$. However, by \Corref{C:sFinerM} we have $\gt_m \subset \gt_s$. Thus $\gt_s = \gt_m$. Likewise, we have $\gt_i \subset \gt_m = \gt_s$; and, by \Propref{P:sIncCont} and the definition of the inductive limit topology on $V'$ we have $\gt_s \subset \gt_i$. Therefore $\gt_s = \gt_m = \gt_i$.
\end{proof}


\section{Borel Field}
In this section our aim is to discuss the \sfield\ on $V'$ generated by the three topologies (strong, weak, and inductive). We will see that under certain conditions the three \sfield s coincide. The standing assumption throughout this section is that $V$ is a \CHS\ with a countable dense subset $Q_o$. On each $V_n' \subset V'$ define the sets $\F_n({\fr(1/k)})$ for all $k$ as:
\begin{center}
$\F_n({\fr(1/k)}) = \Fnk{k}$,
\end{center}
where $Q = Q_o - \{0\}$.

Recall that the local base for the topology of $V_n'$ is given by the sets
\begin{center}
$\U_n(\e) = \Vnball{\e}$,
\end{center}
where $\e > 0$.

\begin{Lemma} \label{L:FnkOpen}
In $V_n'$ we have that $\F_n({\fr(1/k)}) = \U_n({\fr(1/k)})$ for all $k$.
\end{Lemma}

\begin{proof}
Recall that $\norm|v'|_{-n} = \sup_{v \in V_n-\{0\}} |\ip<v', \tfr(v/{\norm|v|_n})>|$.
It is enough to show that for any $v' \in V_n'$ we have $\norm|v'|_{-n} = \sup_{v \in Q} |\ip<v', \tfr(v/{\norm|v|_n})>|$. This is quite easy to see: for any non-zero $v \in V_n$ we have a sequence $\{v_l \}_{l = 1}^{\infty}$ in $Q$ that converges to $v$ (since $Q$ is dense in $V$ and $V$ is dense in $V_n$). Thus $\ip<v',v> = \lim_{l \to \infty} \ip<v', v_l>$. 
\end{proof}

\begin{Prop} \label{P:fCount}
The collection $\{ \F_n({\fr(1/k)}) \}_{k=1}^{\infty}$ forms a local base in $V_n'$. That is, $V_n'$ is first countable.
\end{Prop}

\begin{proof}
Take an open set $U \subset V_n'$ containing $0$. Then $\U_n(\e) \subset U$ for some $\e > 0$. Choose $k$ so that $\fr(1/k) < \e$. Then by \Lemmaref{L:FnkOpen} we have $\F_n({\fr(1/k)}) = \U_n({\fr(1/k)}) \subset \U_n(\e)$.
\end{proof}

Since each $V_n$ is a separable Hilbert space, so is its dual $V_n'$. Let $Q_n'$ be a countable dense subset in $V_n'$.

\begin{Prop} \label{P:sCount}
The collection $\{x' +  \F_n({\fr(1/k)}) \medspace \mid \medspace x' \in Q_n' \, , \, 1 \leq k < \infty \}$ is a basis for $V_n'$. That is, $V_n'$ is second countable.
\end{Prop}

\begin{proof}
Consider an open set $U \subset V_n'$ and an element $v'$ in $U$. By \Propref{P:fCount} there is a $k$ such that $v' + \F_n({\fr(1/k)}) \subset U$. Take $x' \in Q_n'$ such that $\norm|x' -v'|_{-n} < \fr(1/2k)$. 

Observe that $x' + \F_n({\fr(1/2k)}) \subset v' + \F_n({\fr(1/k)})$: Take any $w' \in \F_n({\fr(1/2k)})$ and we have
\begin{align*}
\sup_{v \in Q} |\ip<x' -v' + w', \tfr(v/{\norm|v|_n})>| &\leq \norm|x'-v'|_{-n} + \sup_{v \in Q} |\ip<w',\tfr(v/{\norm|v|_n})>| \\
&< \fr(1/2k) + \fr(1/2k) = \fr(1/k).
\end{align*}
This gives us that $x' -v' + \F_n({\fr(1/2k)}) \subset  \F_n({\fr(1/k)})$ or equivalently
$x' + \F_n({\fr(1/2k)}) \subset v' + \F_n({\fr(1/k)})$. Also $v' \in x' + \F_n({\fr(1/2k)})$ since $\norm|x'-v'|_{-n} < \fr(1/2k)$. 

In summary we have that $v' \in x' + \F_n({\fr(1/2k)}) \subset v' + \F_n({\fr(1/k)}) \subset U$. Therefore the collection $\{x' +  \F_n({\fr(1/k)}) \medspace \mid \medspace x' \in Q_n' \, , \, 1 \leq k < \infty \}$ is a basis for $V_n'$
\end{proof}

\begin{Lemma} \label{L:FnkBorel}
Let $\sigma (\gt_w)$ be the Borel \sfield\ on $V'$ induced by the weak topology. Then $F_n({\fr(1/k)})$ is in $\sigma (\gt_w)$  for all positive integers $k$ and $n$.
\end{Lemma}

\begin{proof}
Observe $\F_n({\fr(1/k)}) = \Vnball{\fr(1/k)} = \Vball{\fr(1/k)}$. (If $v' \in V'$ satisfies $\sup_{v \in Q} \vert \ip<v', \fr(v/{\norm|v|_n})> \vert < \fr(1/k)$, then $v' \in V_n'$.)

Now note that $\F_n({\fr(1/k) })$ can be expressed as
\[
\F_n({\tfr(1/k)}) = \bigcup_{r \in S}  \bigcap_{v \in Q_n'}  \bN({\tfr(v/\norm|v|_n)}; r),
\]
 where $\bN({\tfr(v/\norm|v|_n)};r) = \rweakS[{\tfr(v/\norm|v|_n)}]{r}$ and $S = \{ r \in \Q \medspace ; \medspace 0 < r < \fr(1/k) \}$. Therefore $\F_n({\fr(1/k)})$ can be expressed as the countable intersection of the weakly open sets $\bN({\fr(v/\norm|v|_n)}; r)$. Hence $\F_n({\fr(1/k)})$ is in $\sigma (\gt_w)$. 
\end{proof}

\begin{Thm} \label{T:borel}
Let $V'$ be endowed with a topology $\gt$. If  $\gt$ is finer than $\gt_w$ and the inclusion map $i_n':V_n' \to V'$ is continuous for all $n$, then the \sfield s generated by $\gt$ and $\gt_w$ are equal. \emph{(}i.e. $\sigma (\gt_w) =  \sigma (\gt)$\emph{)} 
\end{Thm}

\begin{proof}
Let $U$ be a set in $\gt$. Then $U_n = i_n'^{-1}(U)$ is open in $V_n'$. By \Propref{P:sCount}, $U_n$ can be expressed as $U_n = \bigcup_{l \in T} x_{n_l}' + F_n({\fr(1/k_l)})$ where $x_{n_l}' \in Q_n'$ and $T$ is countable. Then
\begin{align*}
  U \cap V' &= U \cap \Big( \bigcup_{n =1}^{\infty} V_n' \Big) = \bigcup_{n =1}^{\infty} U \cap V_n' \\
  &= \bigcup_{n =1}^{\infty} U_n = \bigcup_{n =1}^{\infty} \bigcup_{l \in T} x_{n_l}' + \F_n({\tfr(1/k_l)}).
\end{align*}
Thus $U$ can be expressed as a countable union of sets in $\sigma(\gt_w)$. Hence $U$ is in $\sigma(\gt_w)$. Therefore $\sigma(\gt_w) = \sigma(\gt)$.

\end{proof}

\begin{Cor}
The \sfield s generated by the inductive, strong, and weak topologies on $V'$ are equivalent. \emph{(}i.e. $\sigma (\gt_w) = \sigma (\gt_s) = \sigma (\gt_i)$\emph{)}
\end{Cor}

\begin{proof}
We can apply \Thmref{T:borel} since $i_n'$ is continuous with respect to $\gt_i$ and $\gt_s$ and also both $\gt_i$ and $\gt_s$ are finer than $\gt_w$.
\end{proof}

The \sfield\ on $V'$ generated by the weak, strong, or inductive topology is referred to as the \emph{Borel field} on $V'$.


\section{A Word on Nuclear Spaces}

Let $V$ be a \CHS\ associated with an increasing sequence of inner-product norms 
$\{ \norm|\cdot|_{n} ; n \geq 1 \}$. Again let $V_n$ be the completion of $V$ with respect to the norm $\norm|\cdot|_n$. 

\begin{Def}
The \CHS\ $V$ is called a \emph{nuclear space} if for any $n$, there exists $m \geq n$ such that the inclusion map from $V_m$ into $V_n$ is a Hilbert-Schmidt operator (i.e. there is an orthonormal basis $\{ e_k \}_{k=1}^{\infty}$ for $V_m$ such that $\sum_{k=1}^{\infty} \norm|e_k|_n^2 < \infty$).
\end{Def}

\begin{Remark}
Note that a trace class operator is also a Hilbert-Schmidt operator and that the product of two Hilbert-Schmidt operators is a trace class operator. Thus $V$ is a nuclear space if and only if for any $n$, there exists $m \geq n$ such that the inclusion map from $V_m$ into $V_n$ is a trace class operator.
\end{Remark}

\begin{Prop}
Let $V$ be a perfect space. Then $V$ has a countable dense subset \emph{(}i.e. $V$ is separable\emph{)}.
\end{Prop}

\begin{proof}
Recall that $V = \bigcap_{n=1}^{\infty} V_n$. We can divide this into two cases: either each $V_n$ is separable or there exists a $k$ such that $V_k$ is not separable.

In the first scenario, since $V \subset V_1$ and $V_1$ is separable, we can find a countable set $Q_1 \subset V$ such that $Q_1$ is dense in $V$ in the norm $\norm|\cdot|_1$. Likewise, we can find $Q_2 \subset V$ that is dense in $V$ with respect to the norm $\norm|\cdot|_2$. Continuing in this manner, we form $Q_n \subset V$ for all $n \in \{1,2,\ldots\}$. Let $Q = \bigcup_{n=1}^{\infty} Q_n$. We will now show $Q$ is dense in $V$. Let $v \in V$. For each $n$ we can find a $v_n \in Q_n$ such that $\norm|v - v_n|_n < \fr(1/n)$. Then for any $k < n$ we have that
\[
\norm|v-v_n|_k \leq \norm|v-v_n|_n < \fr(1/n)
\]
Therefore, the sequence $\{v_n\}_{n=1}^{\infty}$ will converge to $v$ in the space $V$.

In the second case, without loss of generality we can take $V_1$ to be nonseparable. Using the Axiom of Choice we can find an uncountable set $S_1$ in $V$ of points bounded in the norm $\norm|\cdot|_1$ with the distance between any two points being larger than a positive constant $M$. (That is, for $x,y \in S_1$, we have $\norm|x-y|_1 \geq M$.) Likewise, since 
$V = \bigcup_{n=1}^{\infty} \{ v \in V \, ; \, \norm|v|_2 \leq n \}$, there is an uncountable set $S_2 \subset S_1$, which is bounded in the norm $\norm|\cdot|_2$. Continuing in this manner, for each $n$ we form an uncountable set $S_n \subset S_{n-1}$ such that $S_n$ is bounded in the norm $\norm|\cdot|_n$. Note that for any $x,y \in S_n$, we have that 
\[
\norm|x-y|_n \geq \norm|x-y|_1 \geq M
\]

From each $S_k$ take an arbitrary point $v_k$ and form the set $\{v_k\}_{k=1}^{\infty}$. Note that $\{v_k\}_{k=1}^{\infty}$ is bounded in $V$. However, by construction $\{v_k\}_{k=1}^{\infty}$ cannot contain a Cauchy sequence. Therefore, $V$ cannot be perfect, a contradiction.
\end{proof}

\begin{Prop}
If $V$ is a nuclear space, then $V$ is perfect.
\end{Prop}

\begin{proof}
Let $B$ be a bounded set in $V$. Denote the set $B$ considered as a subset of $V_n$ by $B_n$. Since $B$ is bounded, each $B_n$ is bounded in $V_n$. For $m < n$, let $i_{n,m} :V_n \to V_m$ be the inclusion map. Note that $i_{n,m}(B_n) = B_m$. Since $V$ is a nuclear space the image of the bounded set $B_n$ has compact closure in $V_m$. For $m=1$, taking a sequence of  elements $\{v_k\}_{k=1}^{\infty}$ in $B$, there is a subsequence $\{v_{k_1}\}_{k_1=1}^{\infty}$ that is Cauchy in the norm $\norm|\cdot|_1$. Taking $m = 2$, we can find subsequence $\{v_{k_2}\}_{k_2=1}^{\infty}$ of $\{v_{k_1}\}_{k_1=1}^{\infty}$ that is Cauchy in the norm $\norm|\cdot|_2$. Continuing in this way and forming the diagonal sequence $\{v_{k_j}\}_{j=1}^{\infty}$ we see that $\{v_{k_j}\}_{j=1}^{\infty}$ is Cauchy in every norm $\norm|\cdot|_k$. Thus 
$\{v_{k_j}\}_{j=1}^{\infty}$ is Cauchy in $V$. Since $V$ is complete, this sequence has a limit in $V$. Thus $B$ is limit point compact.
\end{proof}

Combining the last two propositions, we see that all the results proved throughout this article apply to nuclear spaces. Most importantly, for a nuclear space, the strong and inductive topologies on the dual coincide and the \sfield s generated by the inductive, strong, and weak topologies are equal.

\section{Gaussian Measure on the Dual of a Nuclear Space}

 Let $E$ be a real separable Hilbert space with norm $\norm|\cdot|_0$, and  let $A$ be a positive Hilbert-Schmidt operator on $E$. Thus $E$ has an orthonormal basis
$\{e_n\}_{n = 1}^{\infty}$ of eigenvectors of $A$,  with
\[
Ae_n=\lambda_ne_n
\]
 and
\[
\sum_{n\geq0}|\lambda_n|^2<\infty \text{ with each } \lambda_n>0
\]

Using the notation $W=\{0,1,2,...\}$, we have the coordinate map
\[
I: E\mapsto \R^W: 
f\mapsto \bigl(\ip<f,e_n> \bigr)_{n\in W}
\] 

Let
\begin{equation}\label{eq:defFzero}
F_0=I(E)=\Big\{(x_n)_{n\in W}: \sum_{n\in W} x_n^2<\infty\Big\}
\end{equation} 
Now, for each $p\in W$, let
\begin{equation}\label{eq:defFp}
F_p=\{(x_n)_{n\in W}: \sum_{n\in W}\lambda_n^{-2p}x_n^2<\infty\}
\end{equation}
On $F_p$ we have the inner-product $\ip<\cdot,\cdot>_p$
given by
$$\ip<a, b>_p =\sum_{n\in W}\lambda_n^{-2p}a_nb_n$$
This makes $F_p$ a real Hilbert space, unitarily isomorphic to
$L^2(W, \gv_p)$ where $\gv_p$ is the measure on $W$ specified by
$\gv_p(\{ n\})=\lambda_n^{-2p}$. Moreover, we have
\begin{equation}\label{eq:chainFp}
F\stackrel{\rm def}{=}\cap_{p\in W}F_p\subset\cdots F_2\subset
F_1\subset F_0=L^2(W,\gv_0)
\end{equation}
and each inclusion $F_{p+1}\to F_p$ is Hilbert-Schmidt.

Now we pull this back to $E$. First set
\begin{equation}\label{eq:defEpnew}
\Ep=I^{-1}(F_p)=\{x\in E: \sum_{n\geq 0}\lambda_n^{-2p}|\ip<
x,e_n>|^2<\infty\} \end{equation} 
It is readily checked that
\begin{equation}\label{eq:EpBp}
\Ep=A^p(E) 
\end{equation} 
On $\Ep$ we have the pull back
inner-product $\ip<\cdot,\cdot>_p$, which works out to
\begin{equation}\label{eq:Epinnerprod}
\ip<f,g>_p=\ip<A^{-p}f, A^{-p}g>
\end{equation} 
Then we have the chain
\begin{equation}\label{eq:chainEp}
\E \stackrel{\rm def}{=}\cap_{p\in W}\Ep\subset\cdots \Ep[2]\subset
\Ep[1] \subset E,
\end{equation}
with each inclusion $\Ep[p+1]\to \Ep$ being Hilbert-Schmidt.

Equip $\E$   with the topology generated by the norms
$\norm|\cdot|_p$ . Then $\E$ is, by definition, a
nuclear space.
 The vectors $e_n$ all lie in $\E$ and the set of all
 rational-linear combinations of these vectors produces a
 countable dense subspace of $\E$. 
 Since $\E$ is a nuclear space,
 the topological dual $\Estar$ is the union of the duals $\Epstar$. In
 fact, we have:
\begin{equation}\label{eq:chainEpdual}
\Estar=\cup_{p\in W}\Epstar \supset \cdots \Epstar[2] \supset \Epstar[1]\supset E'\simeq E,
\end{equation}
where in the last step we used the usual Hilbert space isomorphism
between $E$ and its dual $E'$.

Going over to the sequence space, $\Epstar$ corresponds to
\begin{equation}\label{eq:defFpprime}
F_{-p}\stackrel{\rm def}{=}
\{(x_n)_{n\in W}: \sum_{n\in W}\lambda_n^{2p}x_n^2<\infty\}
\end{equation}
The element $y\in F_{-p}$ corresponds to the linear functional on
$F_p$ given by
$$x\mapsto \sum_{n\in W}x_ny_n$$
which, by Cauchy-Schwartz, is well-defined and does define an
element of the dual $F'_p$ with norm equal to the square root of $\sum_{n\in W}\lambda_n^{2p}y_n^2<\infty$.

Consider now the product space $\R^W$, along with the
coordinate projection maps 
\[
{\hat X}_j: \R^W\to {\R}:x\mapsto x_j
\]
 for each $j\in W$. Equip $\R^W$
with the product $\sigma$--algebra, i.e. the smallest sigma-algebra
with respect to which each projection map ${\hat X}_j$ is
measurable. A fundamental result in probability measure theory (a
special case of Kolmogorov's theorem, for instance) says that there
is a unique
 probability measure $\gv$ on the product $\sigma$--algebra  such that
 each function ${\hat X}_j$, viewed as a random variable,
  has standard Gaussian distribution. Thus,
 \[
 \int_{\R^W}e^{it{\hat X}_j}\,d\gv=e^{-t^2/2}
 \]
 for $t\in \R$,  and every $j\in W$. The measure $\gv$ is the
product of the
 standard Gaussian measure $e^{-x^2/2}(2\pi)^{-1/2}dx$ on each
 component ${\R}$ of the product space $\R^W$.

Since, for any   $p\geq 1$, we have
\[
\int_{{\R^W}}\sum_{j\in W}\lambda_j^{2p}x_j^2\,d\nu(x)=\sum_{j\in W}\lambda_j^{2p}<\infty,
\]
it follows that
\[
\gv(F_{-p})=1
\]
for all $p\geq 1$. Thus $\gv(F')=1$.

We can, therefore, transfer  the measure $\gv$ back to $\Estar$,
obtaining a probability measure $\gm$ on the sigma--algebra of
subsets of $\Estar$ generated by the maps
\[
{\hat e}_j:\Estar \to \R:
f\mapsto f(e_j),
\]
where $\{e_j\}_{j\in W}$ is the orthonormal basis of $E$ we
started with (note that each $e_j$ lies in $\E=\cap_{p\geq 0}\Ep$).
This is clearly the sigma--algebra generated by the weak topology on
$\Estar$, which is equal to the sigma--algebras
generated by the strong or inductive-limit topologies.

The above discussion gives a simple direct description of the
measure $\gm$. Its existence is also obtainable by applying the
well--known Minlos theorem.

To summarize, are at the starting point of much of
infinite--dimensional distribution theory (white noise analysis):
Given a real, separable Hilbert space $E$ and a positive
Hilbert-Schmidt operator $A$ on $E$, we have constructed a nuclear
space $\E$ and a unique probability  measure $\gu$ on the Borel
sigma--algebra of the dual $\Estar$ such that there is a linear map
\[
E\to L^2(\Estar,\mu):\xi\mapsto {\hat \xi},
\]
satisfying
\[
\int_{\Estar}e^{it{\hat \xi}(x)}\,d\gu(x)=e^{-t^2\norm|\xi|_0^2/2}
\]
for every real $t$ and $\xi\in E$. This  measure $\gm$ is
often called the (\emph{standard\emph{)} Gaussian measure} or the {\it white noise measure} and is the principal measure used white--noise analysis.

\nocite{*} 
\bibliographystyle{siam}
\bibliography{NuclearSpace}

\end{document}